\documentclass{article}
\usepackage{PRIMEarxiv}
\usepackage[utf8]{inputenc} 
\usepackage[T1]{fontenc}    
\usepackage[colorlinks=true,allcolors=black]{hyperref}      
\usepackage{url}            
\usepackage{booktabs}       
\usepackage{amsfonts}       
\usepackage{amsmath}
\usepackage{mathpazo}
\usepackage{caption}
\usepackage{verbatim}
\usepackage{subcaption}
\usepackage{nicefrac}      
\usepackage{microtype}      
\usepackage{fancyhdr}       
\usepackage{graphicx}       
\graphicspath{{media/}}     
\usepackage{cite}
\usepackage{titlesec}
\usepackage{varwidth}
\title{Singular value decomposition based matrix surgery}   

\author{
  Jehan Ghafuri {\normalfont and} Sabah Jassim \\
  School of Computing \\
  The University of Buckingham \\
  Buckingham, UK, MK18 1EG\\
  \texttt{\{1901699, sabah.jassim\}@buckingham.ac.uk} \\
  }

\begin{document}
\maketitle

\begin{abstract}
This paper aims to develop a simple procedure to reduce and control the condition number of random matrices, and investigate the effect on the persistent homology (PH) of point clouds of well- and ill-conditioned matrices. For a square matrix generated randomly using Gaussian/Uniform distribution, the SVD-Surgery procedure works by: (1) computing its singular value decomposition (SVD), (2) replacing the diagonal factor by changing a list of the smaller singular values by a convex linear combination of the entries in the list, and (3) compute the new matrix by reversing the SVD. Applying SVD-Surgery on a matrix often results in having different diagonal factor to those of the input matrix. The spatial distribution of random square matrices are known to be correlated to the distribution of their condition numbers. The persistent homology (PH) investigations, therefore, are focused on comparing the effect of SVD-Surgery on point clouds of large datasets of randomly generated well-conditioned and ill-conditioned matrices, as well as that of the point clouds formed by their inverses. This work is motivated by the desire to stabilise the impact of Deep Learning (DL) training on medical images in terms of the condition numbers of their sets of convolution filters as a mean of reducing overfitting and improving robustness against tolerable amounts of image noise. When applied to convolution filters during training, the SVD-Surgery acts as a spectral regularisation of the DL model without the need for learning extra parameters. We shall demonstrate that for several point clouds of sufficiently large convolution filters our simple strategy preserve filters norm and reduces the norm of its inverse depending on the chosen linear combination parameters. Moreover, our approach showed significant improvements towards the well-conditioning of matrices and stable topological behaviour.
\end{abstract}


\section{Introduction}
Despite the remarkable success and advancements of deep learning (DL) models in computer vision tasks, serious obstacles to deployment of AI in different domains relates to the challenge of developing deep neural networks that are both robust and generalise well beyond the training data \cite{colbrook2022difficulty}. Accurate and stable numerical algorithms play a significant role to computing a robust and reliable computational models \cite{Higham2002AccuracyAlgorithms}. The source of numerical instability in DL models are partially due to the use of a large number of parameters/hyperparameters, and data that suffer from floating-point errors and inaccurate results. In the case of convolutional neural networks (CNN), an obvious contributor to the instability of their large volume of weights is the repeated action of backpropagation algorithm in controlling the growth of gradient descent to fit the model’s performance to the different patches of training samples. This paper is concerned with empirical estimation of CNN training-caused fluctuation in condition numbers of various weight matrices  as a potential source of instability at convolutional layers and the overall model performance. We shall propose a spectral based approach to reduce and control the undesirable fluctuation. 

The \textit{condition number} $\kappa (A)$ of a square $n \times n$ matrix $A$, considered as a linear transformation $\mathbb{R}^{n \times n} \xrightarrow{} \mathbb{R}$, measures the sensitivity of computing its action to perturbations to input data and round-off errors, defined as $sup \|Ax\|/\|x\|$ over the set of nonzero $x$. It depends on how much the calculation of its inverse suffer from underflow (i.e. how much $det(A)$ is significantly different from 0). Stable action of A means that small changes in the input data are expected to lead to small changes in the output data and these changes are bound by reciprocal of condition number. Hence, the higher the condition number of A is the more unstable its action is to small data perturbation and such matrices are said to be ill-conditioned. Indeed, the distribution of condition numbers of a random matrix simply describes the loss in precision, in terms of the number of digits, as well as the speed of convergence due to ill-conditioning when solving linear systems of equations iteratively, \cite{Edelman1989}. Originally, the condition number of a matrix was first introduced by A. Turing in \cite{turing1948rounding}. Afterwards, the condition number of matrices and numerical problems was comprehensively investigated in \cite{rice1966theory, Demmel1987b, higham1995condition}.  The most common efficient and stable way of computing $\kappa (A)$, is by computing the SVD of $A$ and calculating the ratio of $A$’s largest singular value to its smallest non-zero one, \cite{Klema1980TheApplications}.

J. W. Demmel, in \cite{Demmel1987b}, investigated the upper and lower bounds of the probability distribution of condition numbers of random matrices and showed that the sets of ill-posed problems including matrix inversion, eigenproblems, and polynomial zero finding all have a common algebraic and geometric structure. In particular, Demmel showed that in the case of matrix inversion, the further away a matrix is from the set of noninvertible matrices, the smaller is its condition number. Accordingly, the spatial distribution of random matrices, in their domains, are indicators of the distribution of their condition numbers. These results provide a clear evidence of the viability of our approach to exploit the tools of topological data analysis (TDA) to investigate the condition number stability of point clouds of random matrices. In general, TDA can be used to capture information about complex topological and geometric structures of point clouds in metric spaces with or without prior knowledge about the data (see \cite{Chazal2017} for more detail). Since the early 2000s, applied topology has entered a new era exploiting the persistent homology (PH) tool to investigate the global and local shape of high dimensional datasets. Various vectorisation of persistence diagrams (PD), generated by the PH tool, encode information about both local geometry and global topology of cloud of convolution filters of CNN models, \cite{adams2021topology}. Here we shall attempt to determine the impact of SVD-Surgery procedure on PD’s of point clouds of CNN well- and ill-conditioned convolution filters.   

\textbf{Contribution:} We introduce a singular value decomposition based matrix surgery (SVD-Surgery) technique to modify matrix condition numbers that is suitable for stabilising the actions of ill-conditioned convolution filters on point clouds of image datasets. It decomposes square matrices by SVD factorisation, replaces the smaller singular values, and then reconstruct the original matrix with the resulting singular value diagonal matrix. SVD-Surgery preserves the norm of the input matrix while reducing the norm of its inverse. This means that SVD-Surgery make changes to the PH of the inverse matrices point clouds. We expect that PH analysis of point clouds of matrices (and those of their inverses) can provide an informative understanding of stability behaviour of DL models of image analysis.

\section{Background to the Motivating Challenge }
The ultimate motivation for this paper is related to specific requirements that arose in our challenging investigations of how to “train an efficient slim convolutional neural network models capable of learning discriminating features of Ultrasound (or any radiological) images for supporting clinical diagnostic decisions”. In particular, the developed model’s predictions are required to be robust against tolerable data perturbation and less prone to overfitting effects when tested on unseen data. 

In machine learning and deep learning, vanishing or exploding gradient and poor convergence is generally due to an ill-conditioning problem. The most common approaches to overcome ill-conditioning are regularisation, data normalisation, re-parameterisation, standardisation, and random dropouts. When training a Deep CNN with extremely large datasets of "natural" images, the convolution filter weights/entries are randomly initialised the entries of which are changed through an extensive training procedure using many image batches over a number of epochs at the end of each of which the back-propagation procedure updates the filter entries for improved performance. The frequent updates of filters’ entries result in non-negligible to significant of fluctuation and instability of their condition numbers causing sensitivity of the trained CNN models, \cite{Ghafuri2020a,Ghafuri2021}. CNN models sensitivity are manifested by overfitting, reduced robustness against noise and vulnerability to adversarial attacks, \cite{goodfellow2014explaining}.

Transfer Learning is  common approach in developing CNN models for analysis of US (or other radiology) image datasets whereby the pretrained filters and other model weights of an existing CNN model (trained on natural images) are used as an initialising parameters for retraining. Unfortunately, condition number instabilities increase in transfer learning mode when used for small datasets of non-natural images resulting in suboptimal performance and a model suffer from overfitting.  

\section{Related work}
Deep Learning CNN models involve a variety of parameters, the complexity of which are dominated by the entries of sets of convolution filters at various convolution layers as well as those of the fully connected neural network layers. The norm and/or variance of these parameters are the main factors considered in designing initialisation strategies to speedup training optimisation and improve model performance in machine and deep learning tasks. Currently, most popular CNN architectures initialise these weights using zero-mean Gaussian distributions with controlled layer dependent/independent variances. Krizhevsky et al., \cite{Krizhevsky2017}, uses a constant standard deviation of $0.01$ to initialise weights in each layer. Due to exponentially vanishing/growing gradient and for compatibility with activation functions, \textit{Glorot} \cite{glorot2010understanding}, or \textit{He} \cite{he2015delving}, weights are initialised with controllable variances per layer. For Glorot, the initialised variances depend on the number of in/out neurons, while He initialisation of the variances are closely linked to their proposed parameterised rectified activation unit (PReLU) designed to improve model fitting with little overfitting risk. In all these initialisation strategies, no explicit consideration is given to filter’s condition numbers or their stability during training. In these cases, our investigations found that post training almost all convolution filters are highly ill-conditioned and hence adversely affect their use in transfer learning for non-natural images. More recent attempts to control the norm of the network layer were proposed in GradInit \cite{Zhu2021GradInit:Training} and MetaInit \cite{DauphinMetaInit:Initialize}. These methods can accelerate the convergence while improving model performance and stability. However, both approaches require extra trainable parameters and controlling the condition number during training is not guaranteed. 

Recently, many research works investigated issues closely related to our objectives, by imposing orthogonality conditions on trainable DL model weights. These include orthonormal and orthogonal weight initialisation techniques, \cite{xie2017all,mishkin2015all,saxe2013exact}, orthogonal convolution \cite{Wang_2020_CVPR}, orthogonal regularizer \cite{sinha2018neural}, orthogonal deep neural networks \cite{Jia2021}, and orthogonal weight normalisation \cite{huang2018orthogonal}. Recalling that orthogonal/orthonormal matrices are optimally well conditioned, these publications indirectly support our hypothesis on the link between DL overfitting and condition numbers of learnt convolution filters. Although, instability of weight matrices’ condition numbers are not discussed explicitly, these related work fit into the emerging paradigm of spectral regularisation of NN layers weight matrices. For example, J. Wang et al, \cite{Wang_2020_CVPR}, assert that imposing orthogonality on convolutional filters is ideal for overcoming training instability of DCNN models and improved performance. Furthermore, A. Sinha \cite{sinha2018neural}, point out that ill-conditioned learnt weight matrix contributes to neural network’s susceptibility to adversarial attacks. In fact, their orthogonal regularisation aims to keeping the learnt weight matrix’s condition number sufficiently low, and demonstrate its increased adversarial accuracy when tested on  the natural image datasets of MNIST and F-MNIST. S. Li et al, in \cite{Jia2021}, note that existing spectral regularisation schemes, are mostly motivated to improve training for empirical applications, conduct a theoretical analysis of such methods using bounds the concept of Generalisation Error (GE) measures that is defined in terms of the training algorithms and the isometry of the application feature space. They conclude that optimal bound on GE is attained when each weight matrix of a DNN has a spectrum of equal singular values, and call such models OrthDNNs. To overcome the high computation requirements of strict OrthDNNs, they define approximate OrthDNNs by periodically applying their Singular Value Bounding (SVB) scheme of hard regularisation. In general, controlling weights' behaviour during training has proven to accelerate the training process and reduce the likelihood of overfitting the model to the training set e.g. weight standardisation in \cite{qiao2019micro}, weight normalisation/reparameterization \cite{Salimans2016a}, centred weight normalisation \cite{Huang2017CenteredNetworks}, and using Newton's iteration controllable orthogonalization \cite{Huang2020ControllableDNNs}. Most of the above proposed techniques have been developed specifically to deal with trainable DL models for the analysis of natural images and one may assume that these techniques are used frequently during the training after each epoch/batch. However, none of the known state-of-the-arts DL models seem to implicitly incorporate these techniques. In fact, our investigations of these commonly used DL models revealed that the final convolution filters are highly ill-conditioned, \cite{Ghafuri2020a}. 

Our literature review revealed that reconditioning and regularisation have long been used in analytical applications to reduce/control the ill-conditioning computations noted. In the late 1980’s, E. Rothwell and B. Drachman, \cite{Rothwell1989AProblems}, proposed an iterative method to reduce the condition number in ill-conditioned matrix problem that is based on regularising the non-zero singular values of the matrix. At each iteration, each of diagonal entry in the SVD of matrix is appended with a ratio of a regularising parameter to the singular value. This algorithm is not efficient to be used for our motivating challenge. In addition, the change of the norm is dependent on the regularising parameter. 

In recent years, there has been a growing interest in using TDA to analyse point clouds of various types and complexity datasets. For example, significant advances and insight has been made in capturing local and global topological and geometric features in high dimensional datasets using PH tools, that includes conventional methods, \cite{Turkes2022}. TDA has also been deployed to interpret deep learning and CNNs learning parameters at various layers \cite{Ghafuri2020a,gabrielsson2019exposition,magai2022topology}, and integrating topology-based methods in deep learning \cite{HoferDeepSignatures,rieck2018neural,Ebli2020SimplicialNetworks, Hajij2020, Hu2021, Gonzalez-Diaz2022}. We shall use TDA to assess the spatial distribution of point clouds of matrices/filters (and their inverses) before and after SVD-Surgery for well and ill-conditioned random matrices. 

\section{Topological data analysis} \label{sec:TDA}
In this section, we briefly introduce persistent homology preliminaries and describe the point cloud setting of randomly generated matrices to investigate their topological behaviours.

\textbf{Persistent homology of point clouds:} 
Persistent homology is a computational tool of TDA that encapsulates the spatial distribution of point clouds of data records, sampled from metric spaces, by recording the topological features of a gradually triangulated shape by connecting pairs of data points according to an increasing distance/similarity sequence of thresholds. For a point cloud $X$ and a list $\{\alpha_i\}_0^m$ of increasing thresholds, the shape $S(X)$ generated by this TDA process is a sequence $\{S(X)_i\}_0^m$  of simplicial complexes ordered by inclusion. Vietoris-Rips simplicial complex (VR), is the most commonly used approach to construct $S(X)$ due to its simplicity and Ripser,\cite{Bauer2021Ripser:Barcodes} to construct VR. The sequence of distance thresholds is referred to as a \textit{filtration} of $S(X)$. The topological features of $S(X)$ consists of the number of holes or voids of different dimensions, known as \textit{Bettie} numbers, in each constituents of  $\{S(X)_{i} \}_0^m$. For $j \leq 0$, the $j$-th Bettie number $B_j(S(X)_i)$, are obtained, respectively, by counting  $B_0$ = \#(connected components), $B_1$= \#(empty loops of more than 3 edges), $B_2$= \#(3D cavities bounded by more than 4 faces), etc. Note that the $B_j(S_i(X))$ is the set of generators of the $j$-th singular homology of the simplicial complex $S_i(X)$. The TDA analysis of $X$ with respect to a filtration $\{ \alpha_i\}_0^m$, is based on the persistency of each element of $B_j(S(X)_i)$ as $i \rightarrow m$. Here, persistency of each element is defined as the difference between its birth (first appearance) and its death (disappearance). It is customary, to visibly represent the $B_j(S(X)_i)$, by a vertically stacked set of barcodes, one for each element by the horizontal straight line joining its birth to its death. For more detailed rigours descriptions (see \cite{edelsbrunner2000topological,Ghrist2008,otter2017roadmap}). For simplicity, the barcode set and the PD of the $B_j(S_i(X))$  are referred to by $H_j$.

\begin{figure}[h]
     \centering
     \begin{subfigure}[b]{0.19\textwidth}
         \centering
         \includegraphics[width=0.9\textwidth]{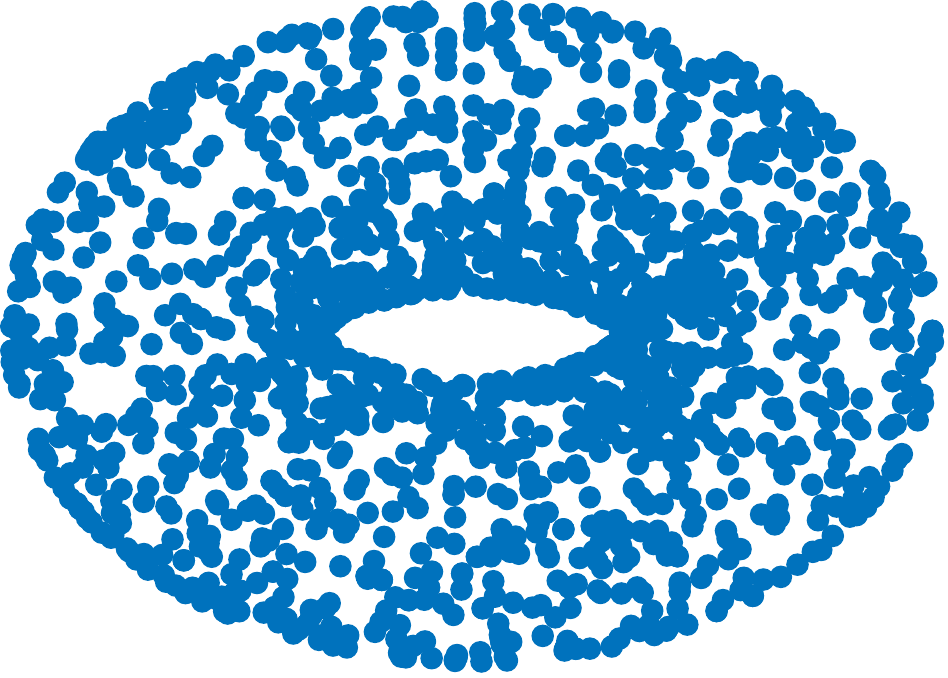}
         \includegraphics[width=\textwidth]{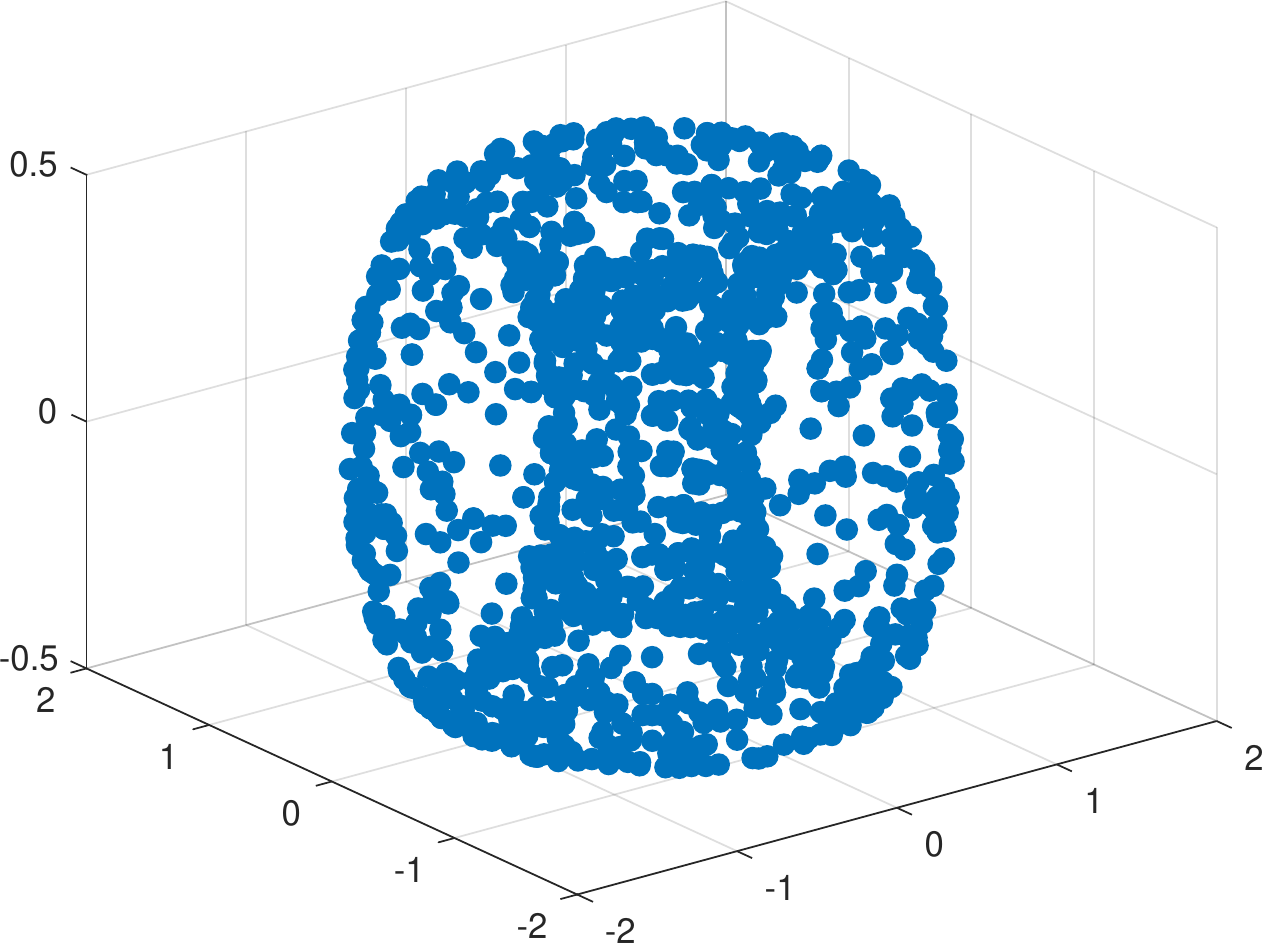}
        \caption{Point cloud}
     \end{subfigure}
     \hspace{0cm}
     \begin{subfigure}[b]{0.19\textwidth}
         \centering
         \includegraphics[width=\textwidth]{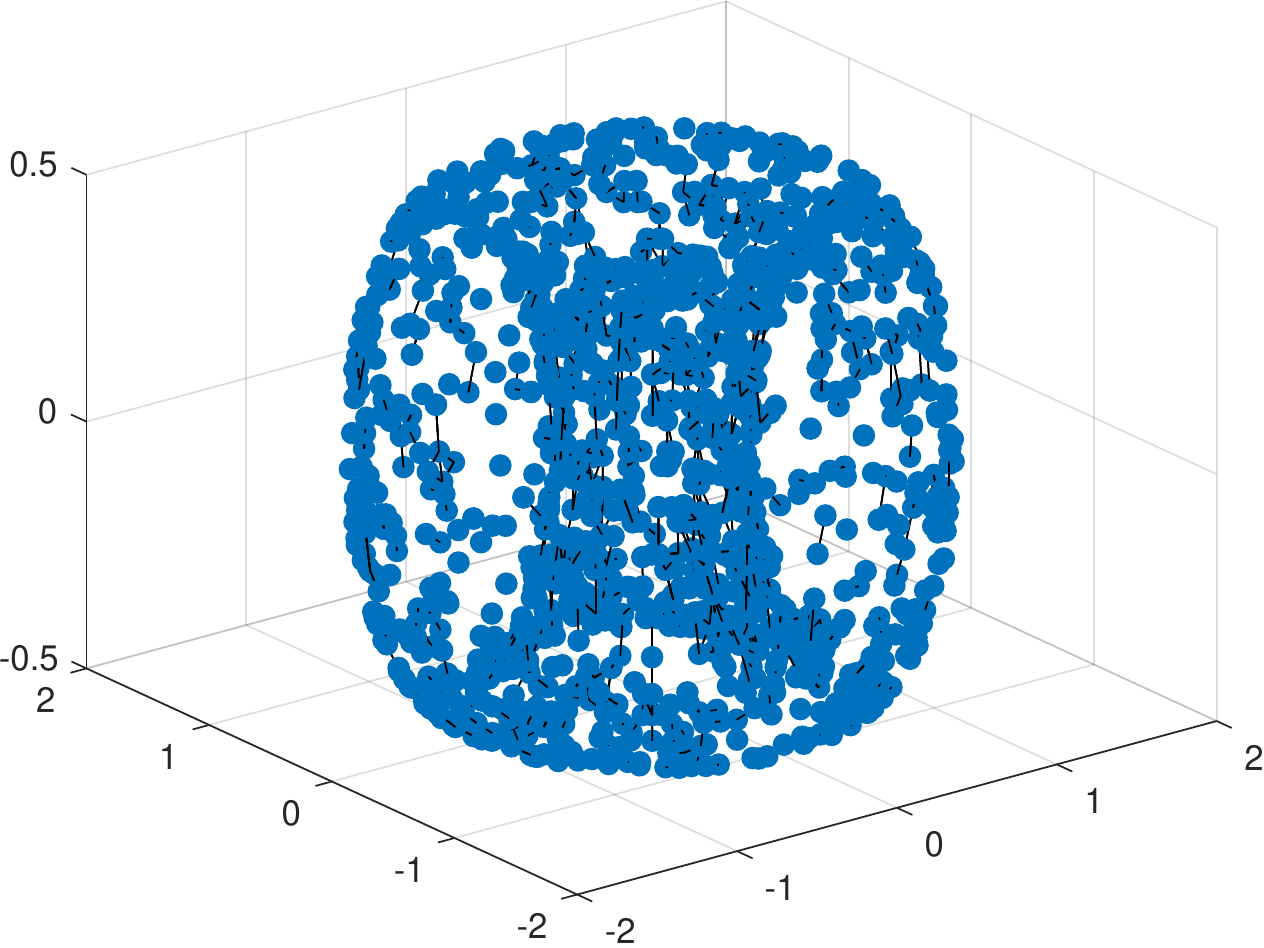}
         \includegraphics[width=\textwidth]{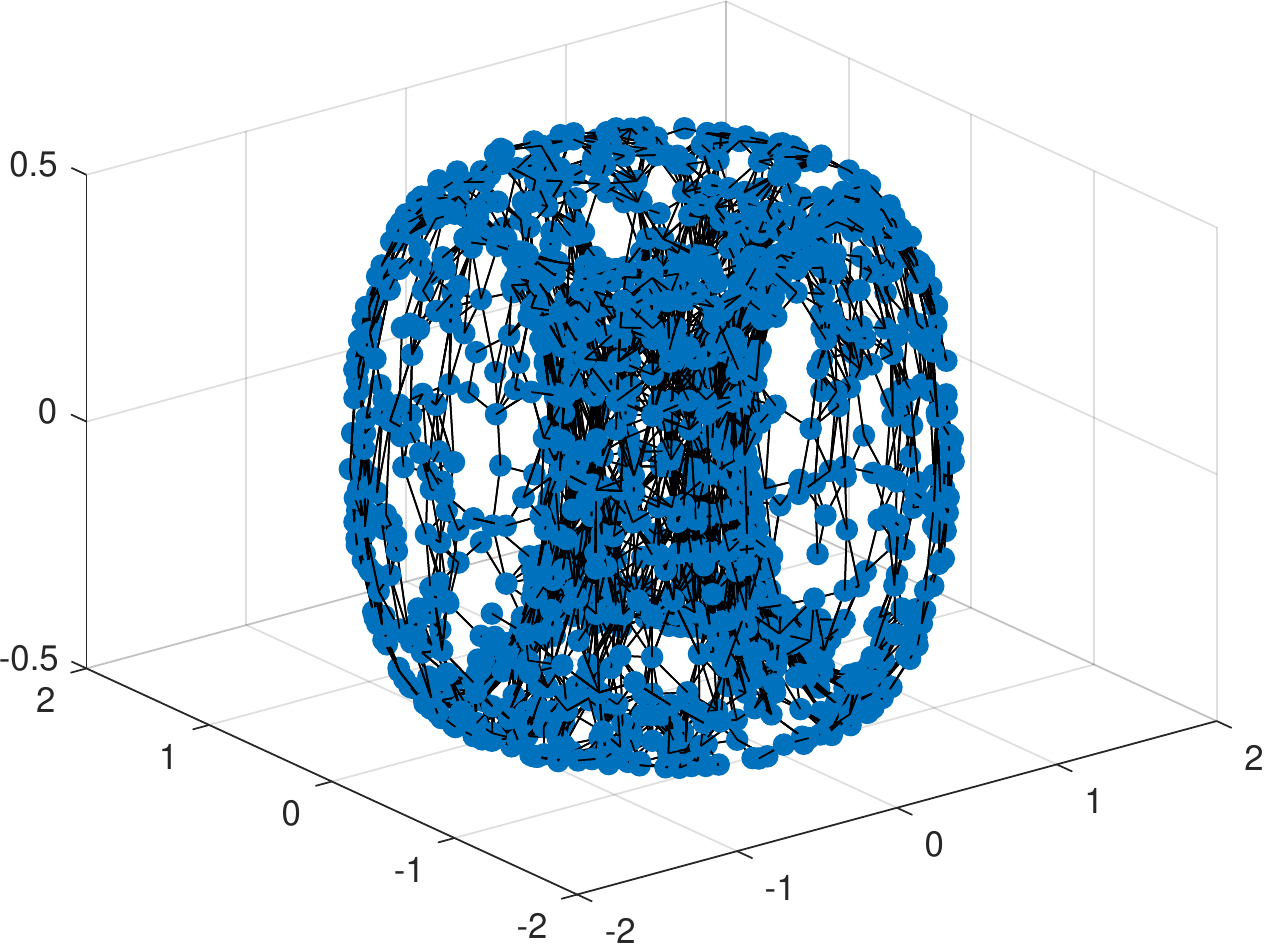}
         \caption{$d \leq 0.1$ \& $0.2$}
     \end{subfigure}
     \hspace{0cm}
     \begin{subfigure}[b]{0.33\textwidth}
         \centering
         \includegraphics[width=\textwidth]{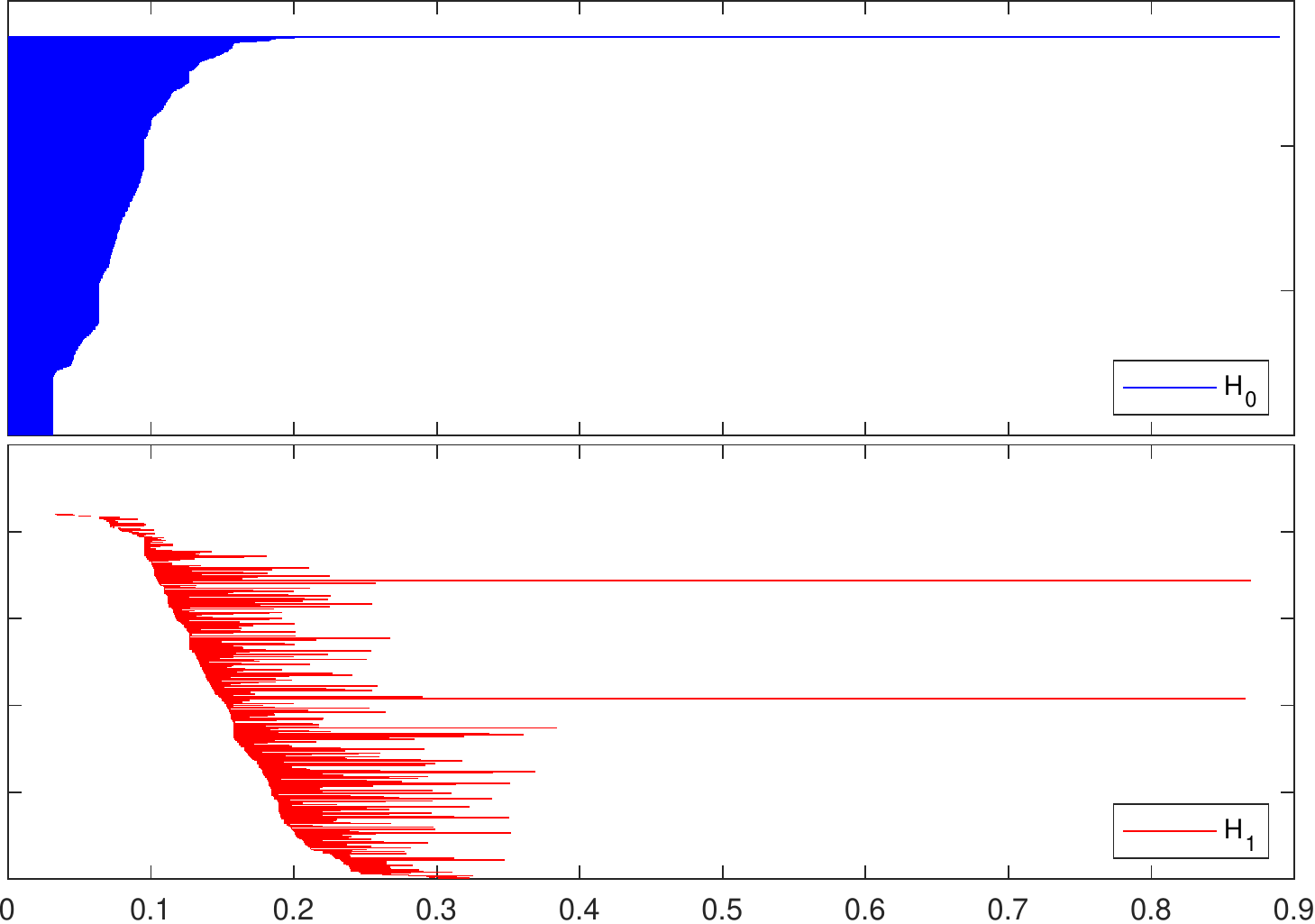}
         \caption{Persistence barcodes}
     \end{subfigure}
     \hspace{0cm}
     \begin{subfigure}[b]{0.26\textwidth}
         \centering
         \includegraphics[width=\textwidth]{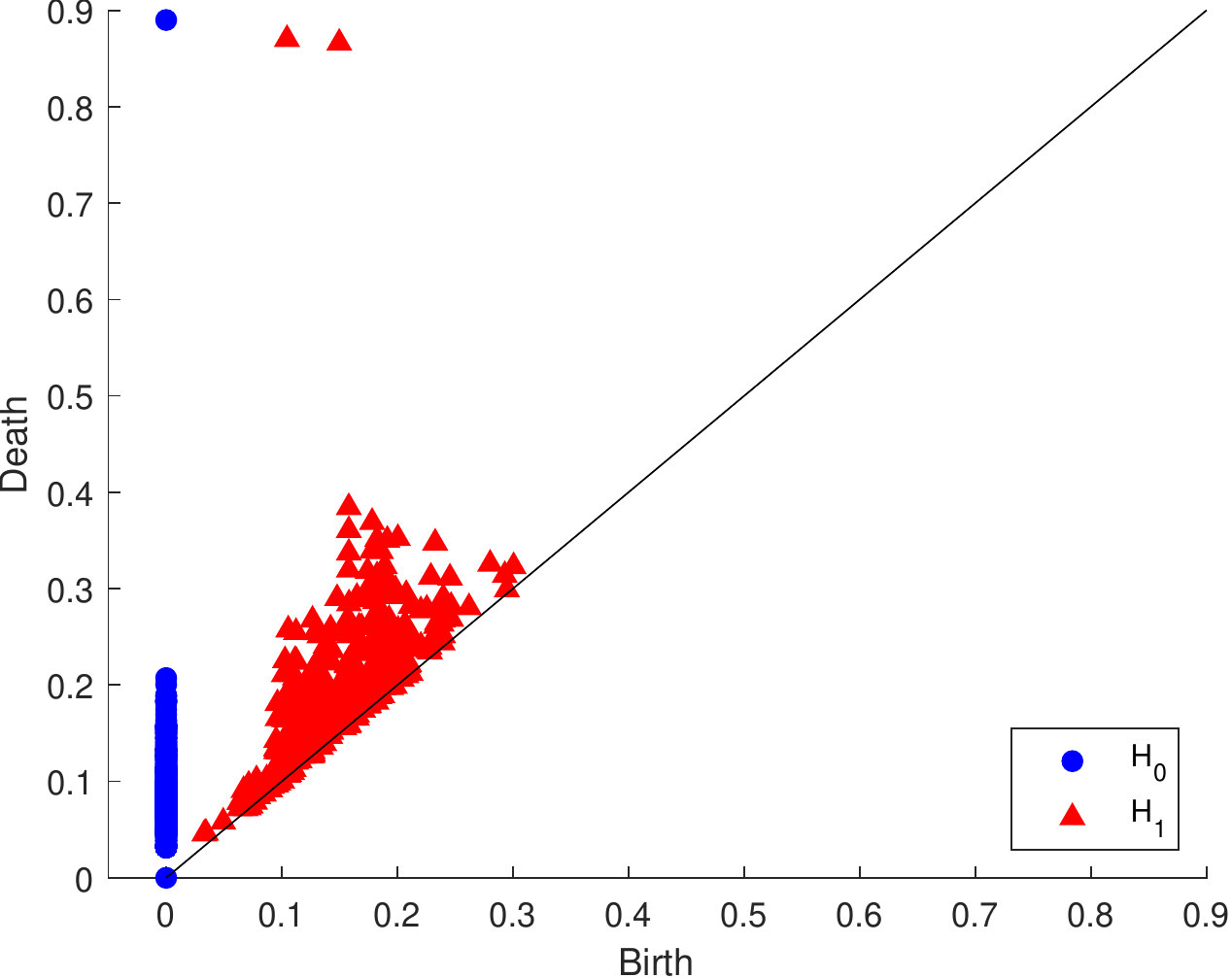}
         \caption{Persistence diagram}
     \end{subfigure}
        \caption{An illustration of a point cloud (a) points from torus, (b) connecting nearby points up to the distances $d$ = $0.1$ \& $0.2$, and their topological representation as persistence barcodes and diagram.} 
       \label{fig:PH_PD_Example}
\end{figure}

Analysis of the resulting PH barcodes of point clouds, in any dimension, is provided by the \textbf{Persistence Diagram} (PD)  formed by a multi-set of points in the first quadrants of the plane $(x=birth$, $y=death)$ above or on the line $y=x$. Each marked point in the PD corresponds to a generator of the persistent homology group of the given dimension, and is represented by a pair of coordinates $(birth,death)$. To illustrate these visual representations of PH information, we created a point cloud of $1500$ points sampled randomly on the surface of the Torus:
\begin{equation*}
    T=\{(x,y,z) \in \mathbb{R}^3 : (\sqrt{x^2+y^2}-a)^2 +z^2=b^2 \}
\end{equation*}
Figure \ref{fig:PH_PD_Example}, below displays this point cloud, together with the barcodes and PD representation of its PH in both dimensions. The two long $1-dim$ persisting barcodes represent the two empty discs whose cartesian product generates the torus. The persistency lengths of these two holes depend on the radii $(a-b, b)$ of the generating circles. In this case, $a=2b$. The persistency lengths of the set of shorter barcodes, are inversely related to the point cloud size.  Noisy sampling will only have an effect on the shorter barcodes.

\begin{figure}[h]
     \centering
     \begin{subfigure}[b]{0.24\textwidth}
         \centering
         \includegraphics[width=\textwidth]{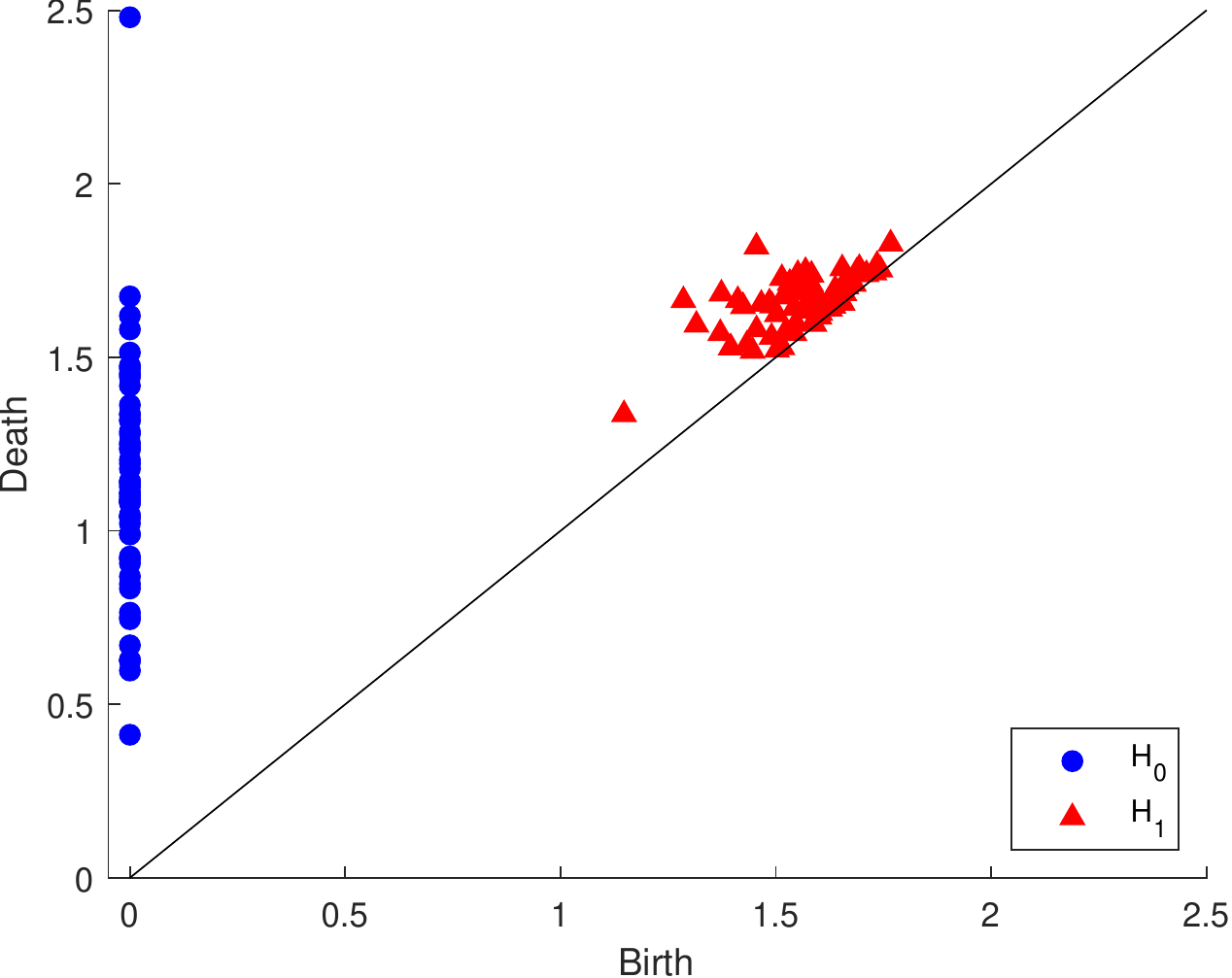}
         \caption{Well-conditioned $\mathcal{X}$}
     \end{subfigure}
      \hspace{0cm}
     \begin{subfigure}[b]{0.24\textwidth}
         \centering
         \includegraphics[width=\textwidth]{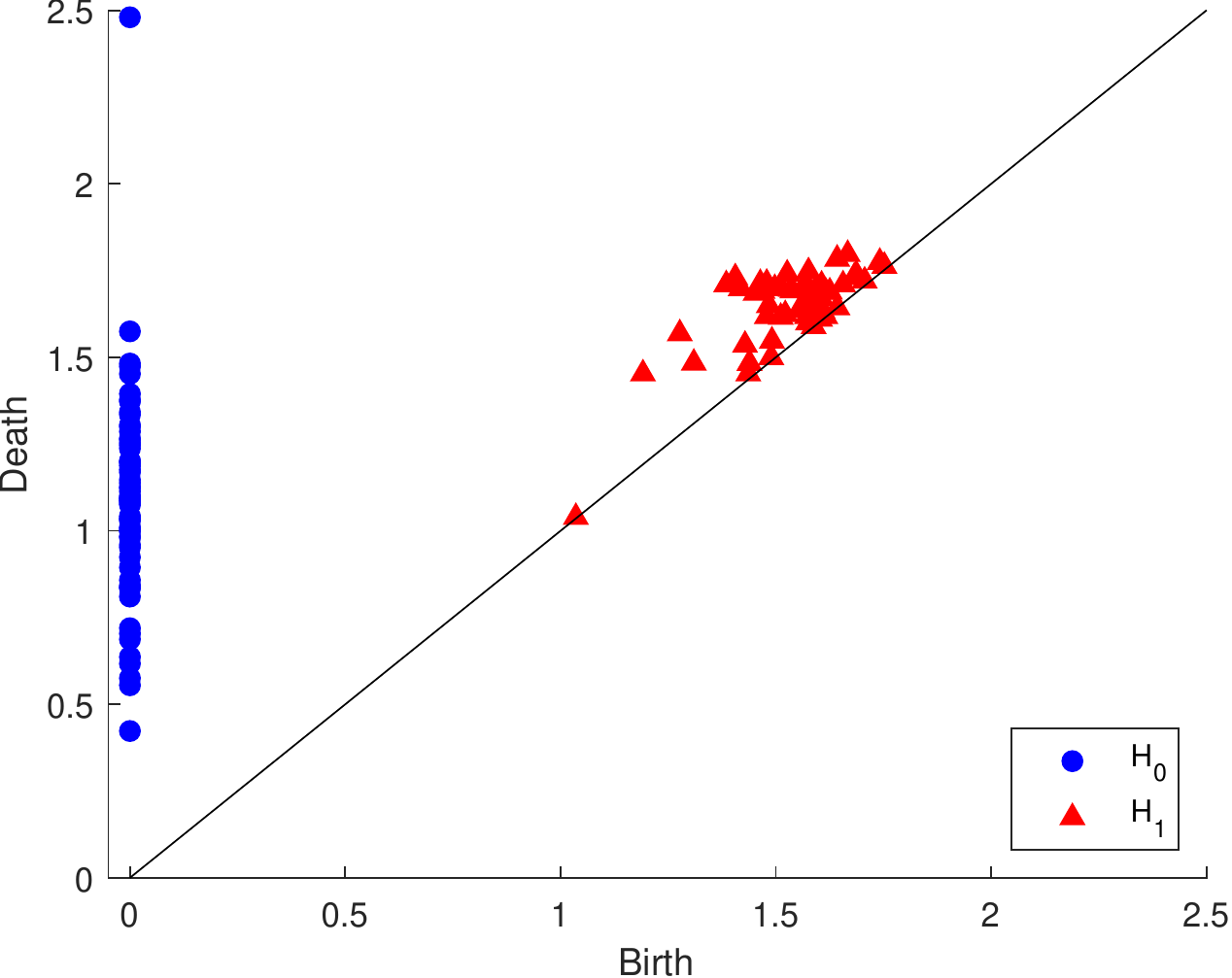}
         \caption{Well-conditioned $\mathcal{X}^{-1}$}
     \end{subfigure}
     \hspace{0cm}
     \begin{subfigure}[b]{0.24\textwidth}
         \centering
         \includegraphics[width=\textwidth]{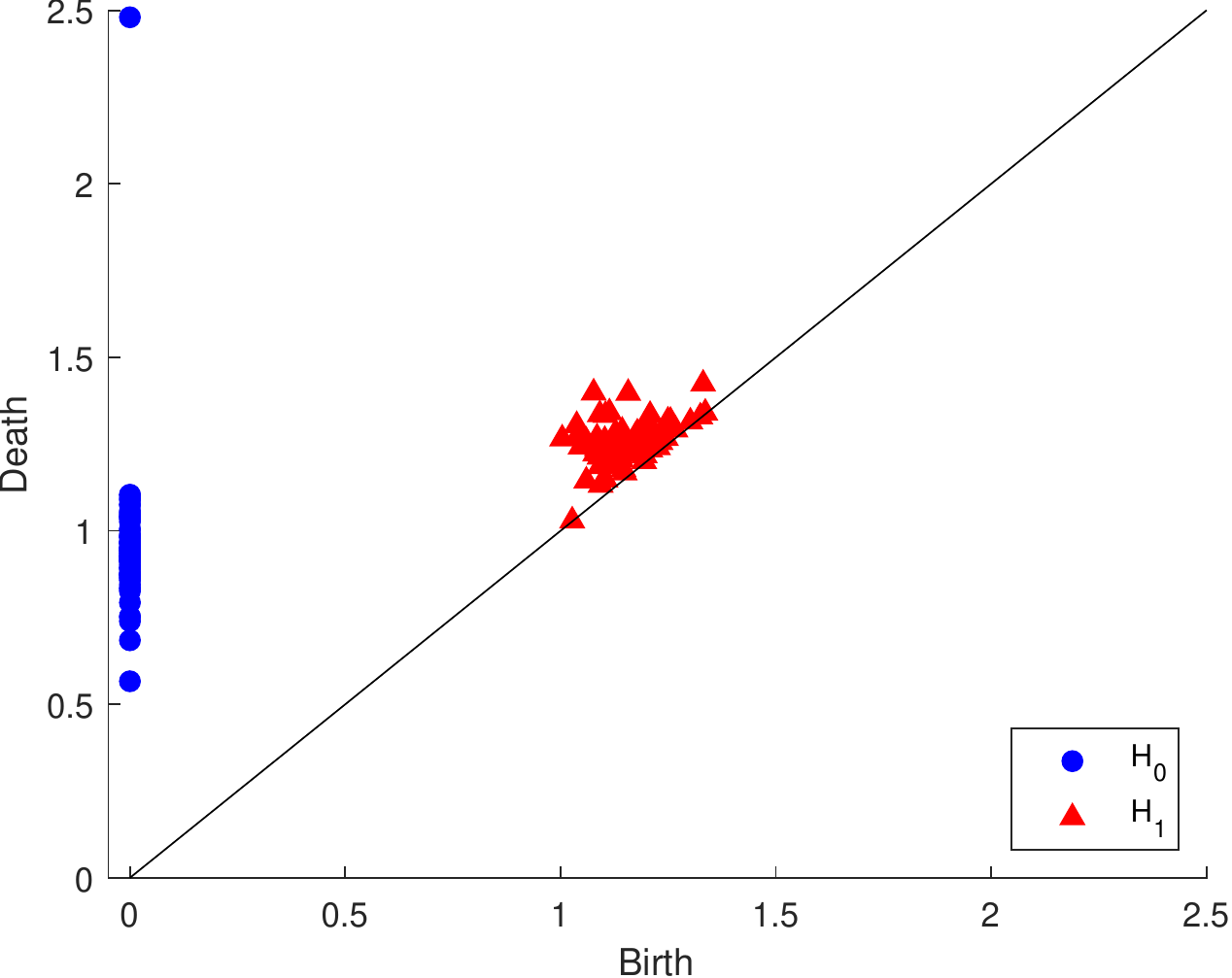}
         \caption{Ill-conditioned $\mathcal{X}$}
     \end{subfigure}
     \hspace{0cm}
     \begin{subfigure}[b]{0.24\textwidth}
         \centering
         \includegraphics[width=\textwidth]{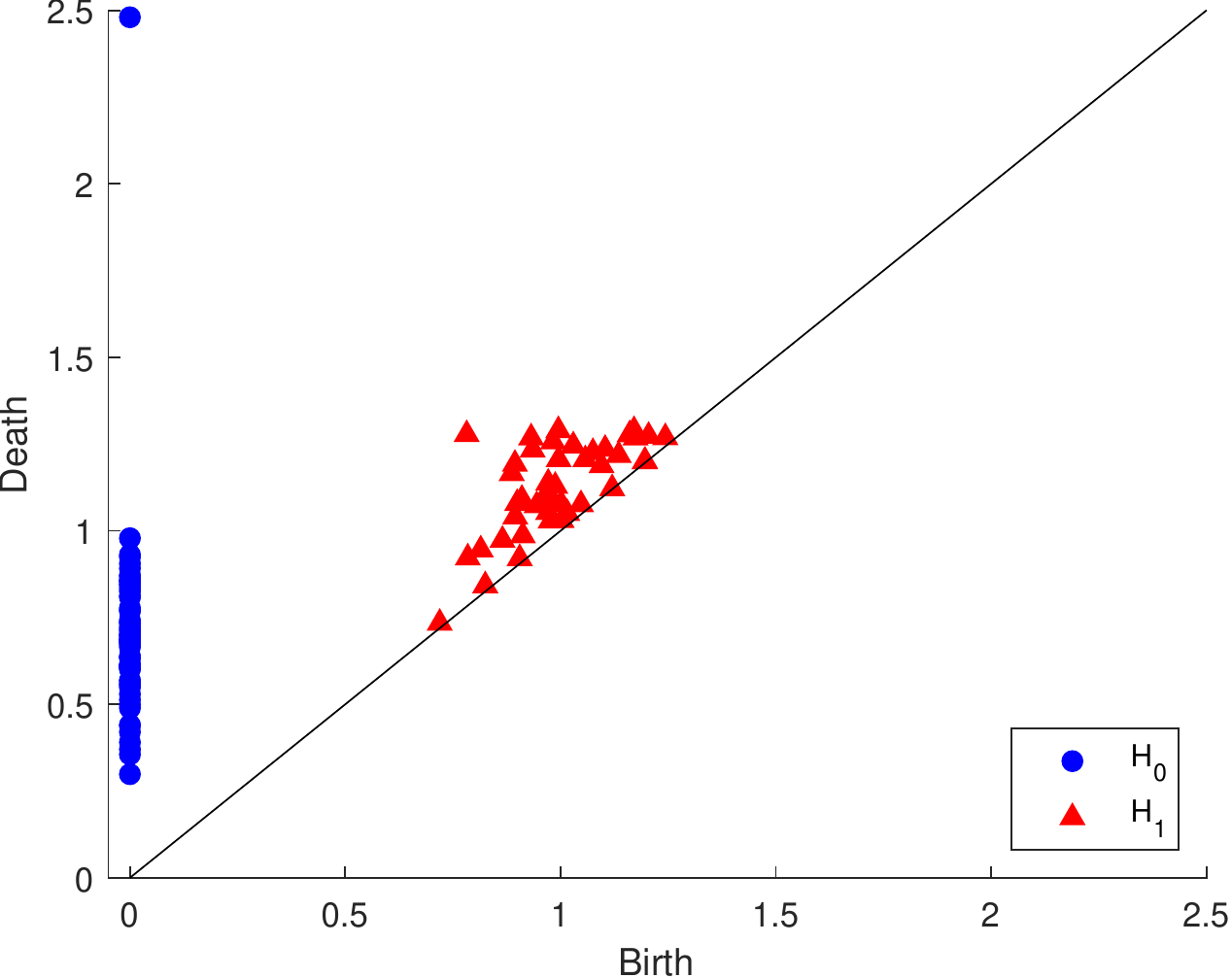}
         \caption{Ill-conditioned $\mathcal{X}^{-1}$}
     \end{subfigure}
        \caption{Persistence diagram of point clouds of the various point clouds representing: (a) well-conditioned matrices, (b) inverse of matrices from (a), (c) ill-conditioned matrices, (d) inverse of matrices from (c).}
        \label{fig:PD well/ill-conditioned filters}
\end{figure}

Our interest in linking PH investigation to our proposed matrix surgery, stems from the significant differences between the topological behaviour (visualised by PDs) of well and ill-conditioned point clouds of filters. Figure \ref{fig:PD well/ill-conditioned filters}, demonstrates this for two randomly generated point clouds of 64 3x3 matrices, and their inverses, according to their lowest and highest condition numbers. 

We note, in both dimensions, the differences between the PDs of the well-conditioned matrices and the ill conditioned ones, as well as between the PD’s of their respective inverses. However, when we examine PDs of the original matrices and their inverse point clouds we find that in dim 0 there is little change in the spatial distributions for the well-conditioned ones compared to that for the ill-conditioned ones. Our proposed matrix surgery will aims to control the differences between the PD’s of the output matrices and that of their inverse point clouds.

\section{Matrix surgery}

In this section, we describe the research framework to perform matrix surgery that aim to reduce and control the condition number of matrices. Suppose matrix $A \in \mathbb{R}^{m \times n}$ is non-singular and based on random Gaussian or Uniform distribution. The condition number of A is defined as: 
\begin{equation*} 
    \kappa (A)=\|A\| \|A^{-1}\|= \frac{\sigma_{1}}{\sigma_{n}}
\end{equation*}
Where $\| . \|$ is the norm of the matrix and we focus on Euclidean norm ($L_2$-norm); $\sigma_{1}$ and $\sigma_{n}$ are the largest and smallest singular values of A, respectively. A matrix is said to be ill-conditioned if any small change in the input results in big changes in the output, and it is said to be well-conditioned if any small change in the input results in a relatively small change in the output. Alternatively, a matrix with a low condition number (close to one) is said to be well-conditioned, while a matrix with a high condition number is said to be ill-conditioned and the ideal condition number of an orthogonal matrix is one. Next , we describe our simple approach of modifying singular value matrix based SVD since the condition number is defined by the largest and smallest singular values. We recall that the Singular Value decomposition of a square matrix $A \in \mathbb{R}^{n \times n}$ is defined by:
\begin{equation*} \label{eq:SVD}
    A = U \Sigma V^T
\end{equation*}
Where $U \in \mathbb{R}^{m \times m}$ and $V \in \mathbb{R}^{n \times n}$ are left and right orthogonal singular vectors (unitary matrices); diagonal matrix $\Sigma=diag(\sigma_1,...,\sigma_n) \in \mathbb{R}^{n \times n}$ are singular values where $\Sigma=\sigma_1 \geq \sigma_2 \geq ... \geq \sigma_n \geq 0$. The SVD-Surgery, described below, is equally applicable to rectangular matrices.

\subsection{SVD based Surgery}

In the wide context, SVD-Surgery refers to the process of transforming matrices to improve their conditioning numbers. In particular, it is targeting matrices that are far from having orthogonality/orthonormality characteristics to replace them with matrices of improved well-conditioned matrix by deploying their left and right orthogonal singular vectors along with the new singular value diagonal matrix. SVD-Surgery can be realised in a variety of ways according to the expected properties of the output matrices to fit the use case. Given any matrix $A$, an SVD-Surgery on $A$ outputs a new matrix of the same size as follows:
\begin{center}
\fbox{\parbox{5.5in}{\centering
\begin{enumerate}
    \item Compute its SVD decomposition,
    \item From the diagonal matrix factor $\Sigma$ construct another diagonal matrix $\tilde{\Sigma}$ by replacing the small singular value(s) while keeping their descendant order
    \begin{equation*} 
    \tilde{\Sigma}= 
    \begin{bmatrix}
          \sigma_1 & 0 & ... & 0\\
          0 & \tilde{\sigma}_2 & ... & 0\\
          \vdots  & \vdots  & \ddots & \vdots  \\
          0 & 0 & ... & \tilde{\sigma}_n
    \end{bmatrix}
    \end{equation*}
    where the updated singular value $\tilde{\sigma}_i$'s are selected to maintain low condition number while the new diagonal entries remain monotonically decreasing, and
    \item Reconstruct the output matrix $\tilde{A}$ as follows:
    \begin{equation*} 
    \tilde{A}= U \tilde{\Sigma} V^T
    \end{equation*}
\end{enumerate}}}
\end{center}
Changes to the singular values amount to rescaling the effect of the matrix action along the left and right orthogonal vectors of $U$ and $V$, and the monotonicity requirement ensures a reasonable control on the various rescaling. While the orthogonal regularisation scheme of \cite{Wang_2020_CVPR} and the SVB scheme of \cite{Jia2021} do reduce condition numbers when applied for improved control of overfitting of DL models trained on natural images, the monotonicity property is not satisfied by these schemes. Moreover, their success cannot be guaranteed for application of DL training of US image datasets. The SVB is much stricter than the SVD-surgery for controlling the condition numbers but no analysis is conducted on the norm of these matrices or their inverses. Our intended SVD-Surgery is designed specifically for use in the motivating application and is aimed to reduce extremely high condition number values and preserve the norm of original matrices. Replacing all singular values with the largest singular value will result in producing an orthogonal matrix with a condition number to be equal to one, but this amounts to ignoring/reducing the effect of significant variations in the training data along some of the singular vectors (i.e. less effective learning). Here, we follow a less drastic strategy in changing singular values:
\begin{center}
\fbox{\parbox{5.5in}{\centering
\vspace{0.08in}
Select a diagonal position $j$, $1<j<n$, let $\tilde{\sigma}$ be a convex linear combination:
    \begin{equation*}
    \tilde{\sigma}_k = \sum_{k=j-1}^{n-1} \alpha_k \sigma_k, \quad where  \quad \alpha_i \geq 0, \quad \sum \alpha_i =1
    \end{equation*}
and set $\tilde{\sigma}_k$  for each  $ j \leq k \leq n.$}}
\vspace{0.1in}
\end{center}
The choice of $j$, and the linear combination parameters can be made application-dependent and possibly determined empirically. 

In the extreme, this strategy include the possibility of setting $\tilde{\sigma}=\sigma_j$. This strategy is rather timid in comparison to the orthogonal regularisation strategies, in that it preserves monotonicity of the singular values. In relation to our motivating application, the parameter choices would be layer dependent, but the linear combination parameters should not result in making significant rescaling of training dataset features along the singular vectors. SVD surgery can be applied to inverse matrices, however, the same replacement strategy and reconstruction may not lead to significant reduction in condition number.
\paragraph{Example} 
Suppose $B$ is a square matrix with $n=3$ drawn from the normal distribution with mean $\mu=0$ and standard deviation $\sigma=0.01$ as follow:
\begin{equation*} \label{matrix:B}
    B= 
    \begin{bmatrix}
    -0.0196	&	0.0291	&	-0.0106\\
    -0.0020	&	0.0083	&	-0.0047\\
    -0.0121	&	0.0138	&	-0.0027
     \end{bmatrix}
\end{equation*}
$B$ can be expressed and decomposed in terms of SVD:

$     U= 
    \begin{bmatrix}
     -0.8745	&	-0.1286	&	-0.4676\\
     -0.2124	&	-0.7652	&	0.6077\\
     -0.4359	&	0.6308	&	0.6419
     \end{bmatrix}
     \Sigma= 
    \begin{bmatrix}
     0.0419	&	0.0000	&	0.0000\\
     0.0000	&	0.0050	&	0.0000\\
     0.0000	&	0.0000	&	0.0005
     \end{bmatrix}
     V= 
    \begin{bmatrix}
     0.5452	&	-0.7156	&	0.4367\\
     -0.7926	&	-0.2703	&	0.5466\\
     0.2731	&	0.6441	&	0.7145
     \end{bmatrix}$
     
$U$ and $V$ are right and left orthogonal singular vectors. Singular values of $B$ are $\Sigma=diag(\sigma_1,\sigma_2,\sigma_3)$ and it is possible to modify and reconstruct $\tilde{B}_1$, $\tilde{B}_2$, and $\tilde{B}_3$ by replacing one and/or two singular values s.t. $\tilde{\Sigma}_1 = diag(\sigma_1,\sigma_2,\sigma_2)$, $\tilde{\Sigma}_2 = diag(\sigma_1,\tilde{\sigma}_2,\tilde{\sigma}_3)$ and $\tilde{\Sigma}_3 = diag(\sigma_1,\sigma_1,\sigma_1)$, respectively. New singular values in $\tilde{\Sigma}_2$ are convex linear combinations s.t. $\tilde{\sigma}_2=2\sigma_1 /3 + \sigma_2 /3$ and $\tilde{\sigma}_3=\tilde{\sigma}_2$ .  
\begin{equation*} \label{matrix:B1 and B2}
    \tilde{B}_1= 
    \begin{bmatrix}
     -0.0205	&	0.0279	&	-0.0121\\
     -0.0008	&	0.0098	&	-0.0027\\
     -0.0108	&	0.0154	&	-0.0007
     \end{bmatrix}
     \tilde{B}_2= 
    \begin{bmatrix}
     -0.0237	&	0.0218	&	-0.0237 \\
      0.0219	&	0.0248	&	-0.0044 \\
     -0.0156	&	0.0204	&	0.0235
     \end{bmatrix}
     \tilde{B}_3= 
    \begin{bmatrix}
    -0.0247	&	0.0220	&	-0.0275 \\
     0.0292	&	0.0296	&	-0.0049 \\
    -0.0171	&	0.0220	&	0.0312
     \end{bmatrix}
\end{equation*}     

After reconstruction, the condition number of $\tilde{B}_1$, $\tilde{B}_2$, and $\tilde{B}_3$ are significantly lower compared to the original matrix as shown in table \ref{tab:Norm and condN}, by using Euclidean norm.

\begin{table}[h]
 \caption{Euclidean norm and condition number of before and after matrix surgery.}
  \centering
  \begin{tabular}{llll}
    \toprule
&	$\|A\|$	&	$\|A^{-1}\|$	&	$\kappa(A)$	\\
	\midrule
$B$	&	0.041883482	&	2034.368572	&	85.20644044	\\
$\tilde{B}_1$	&	0.041883482	&	199.5721482	&	8.358776572	\\
$\tilde{B}_2$	&	0.041883482	&	30.36464182	&	1.271776943	\\
$\tilde{B}_3$	&	0.041883482	&	23.87576058	&	1	\\
    \bottomrule
  \end{tabular}
  \label{tab:Norm and condN}
\end{table}  

\subsection{Effects of SVD-Surgery on large datasets of convolution filters and their inverses}
\begin{figure}[h]
     \centering
     \begin{subfigure}[b]{0.24\textwidth}
         \centering
         \includegraphics[width=\textwidth]{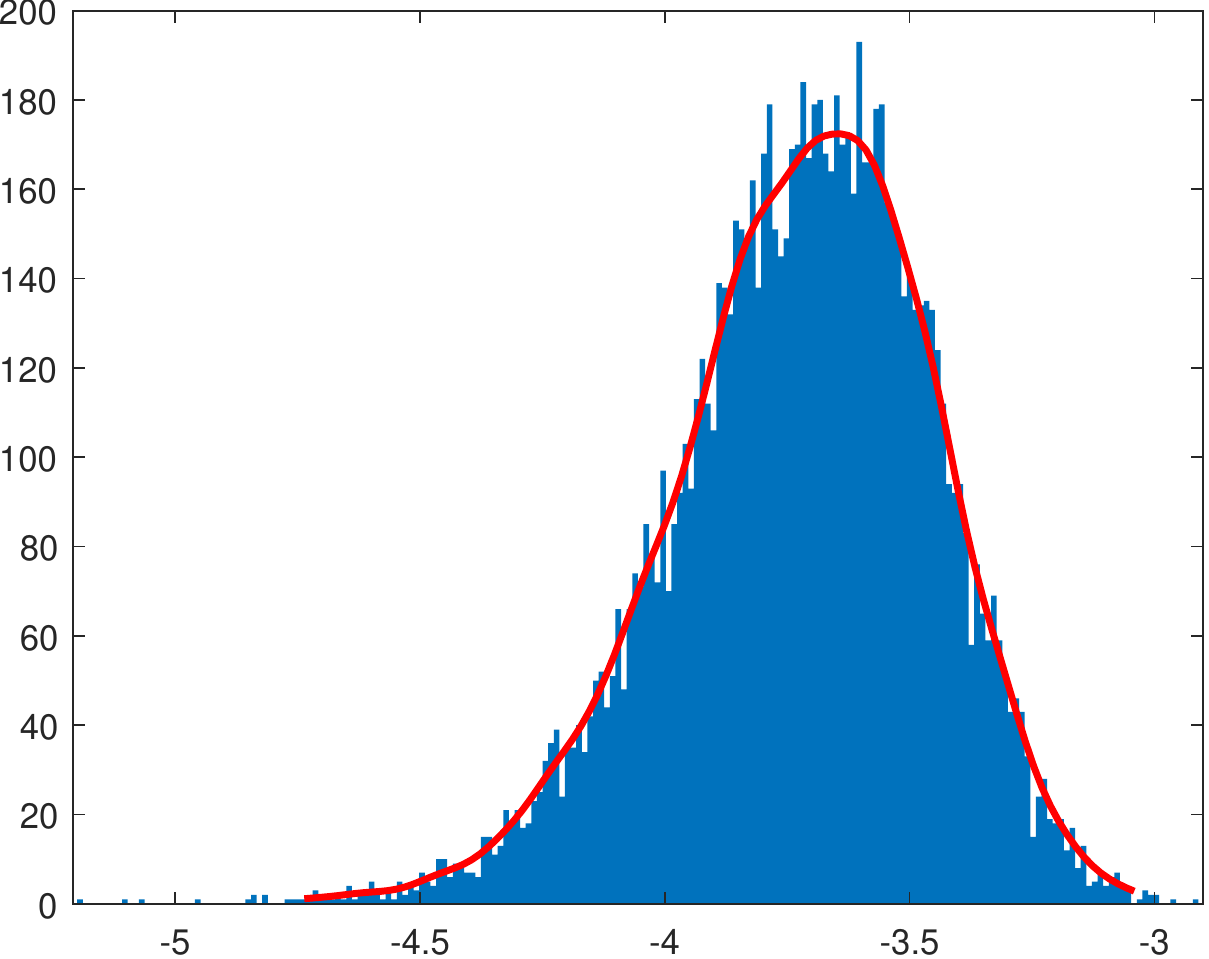}
         \caption{Pre-Surgery $\|A\|$}
         \label{subfig:Pre-Surgery norm}
     \end{subfigure}
     \hspace{0.2cm}
     \begin{subfigure}[b]{0.24\textwidth}
         \centering
         \includegraphics[width=\textwidth]{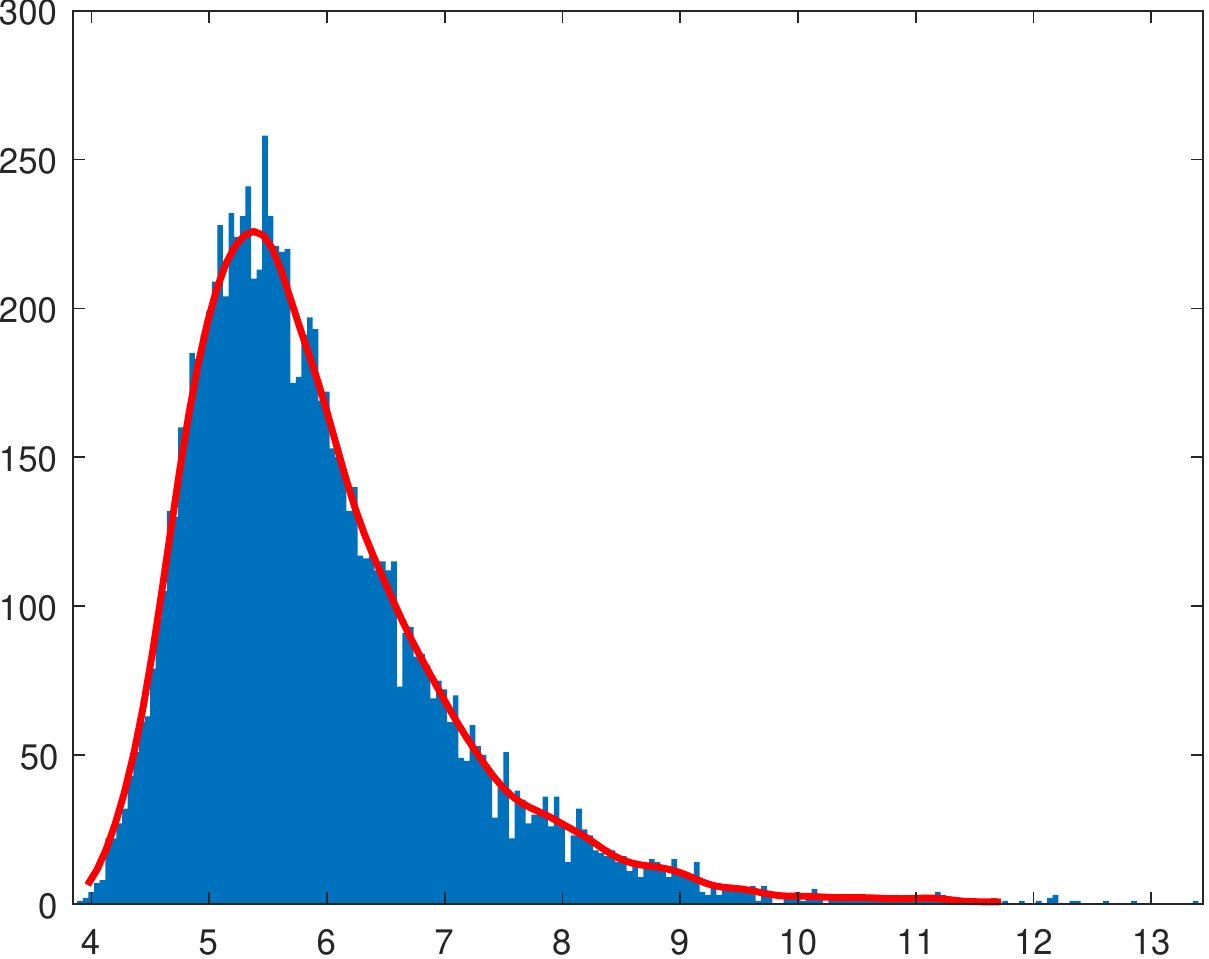}
         \caption{Pre-Surgery $\|A^{-1}\|$}
         \label{subfig:Pre-Surgery normi}
     \end{subfigure}
     \hspace{0.2cm}
     \begin{subfigure}[b]{0.24\textwidth}
         \centering
         \includegraphics[width=\textwidth]{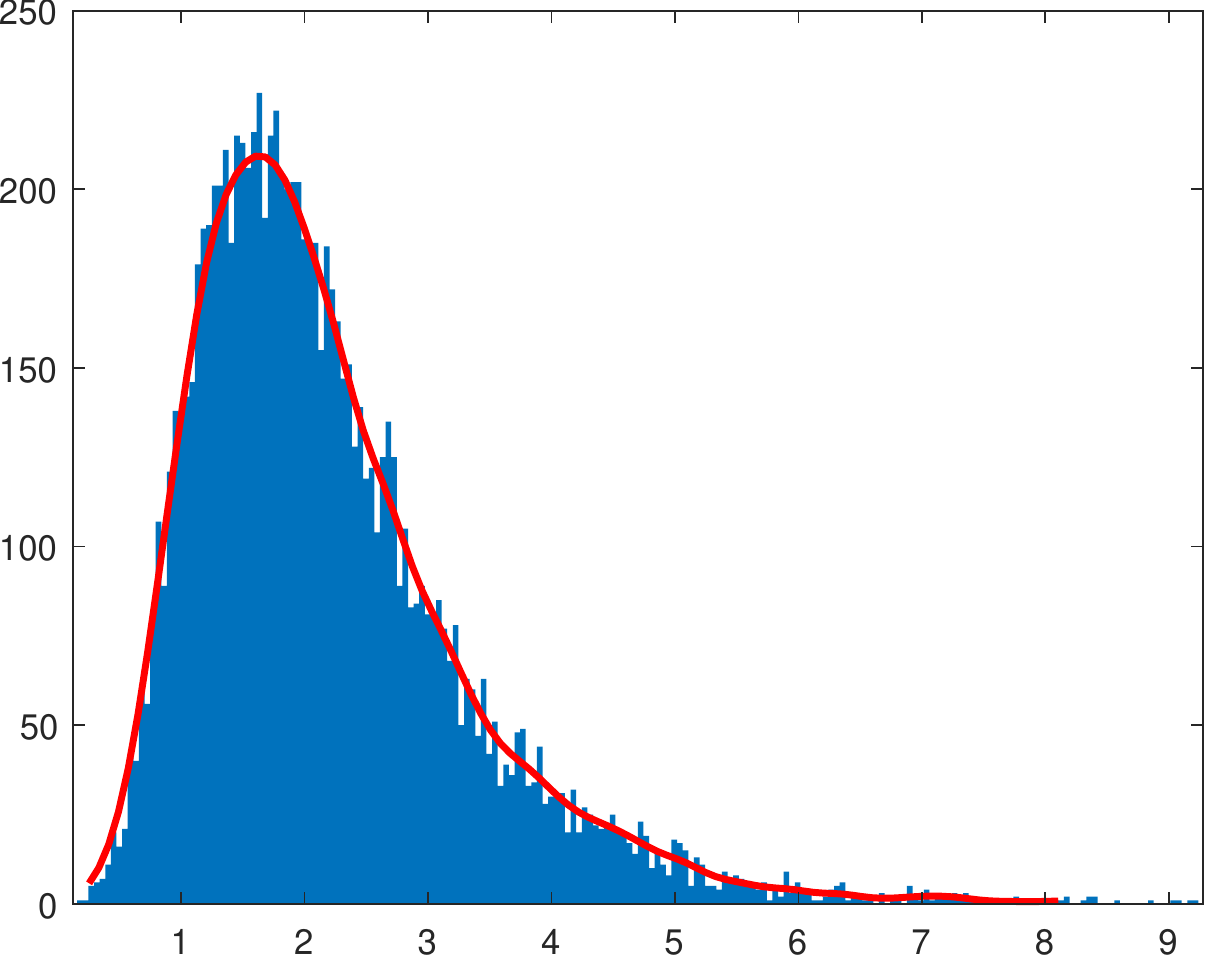}
         \caption{Pre-Surgery $\kappa (A)$}
         \label{subfig:Pre-Surgery CondN}
     \end{subfigure}
     \hspace{0.2cm}
     \begin{subfigure}[b]{0.24\textwidth}
         \centering
         \includegraphics[width=\textwidth]{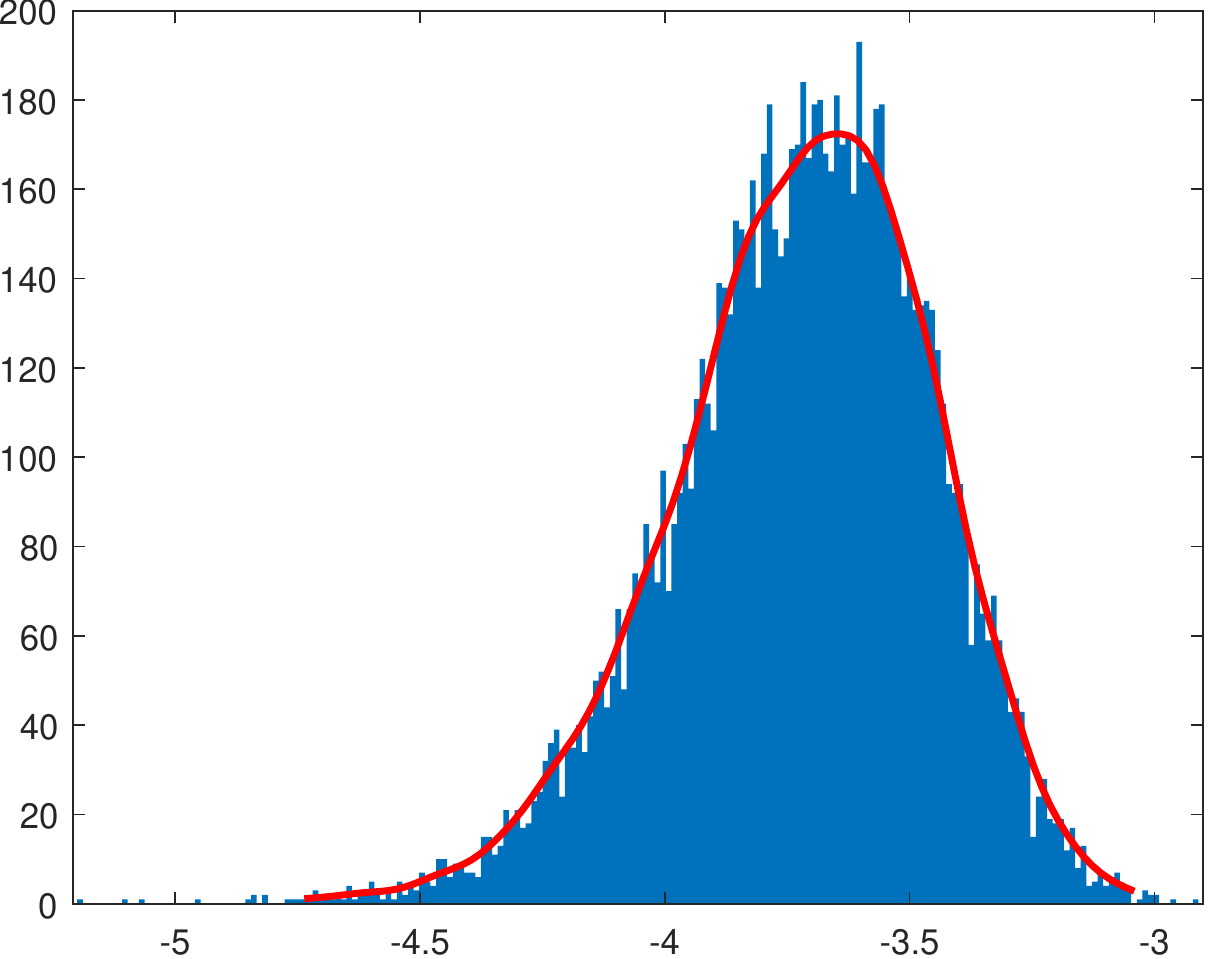}
         \caption{Post-Surgery $\|\tilde{A}\|$}
         \label{subfig:Post-Surgery norm}
     \end{subfigure}
     \hspace{0.2cm}
     \begin{subfigure}[b]{0.24\textwidth}
         \centering
         \includegraphics[width=\textwidth]{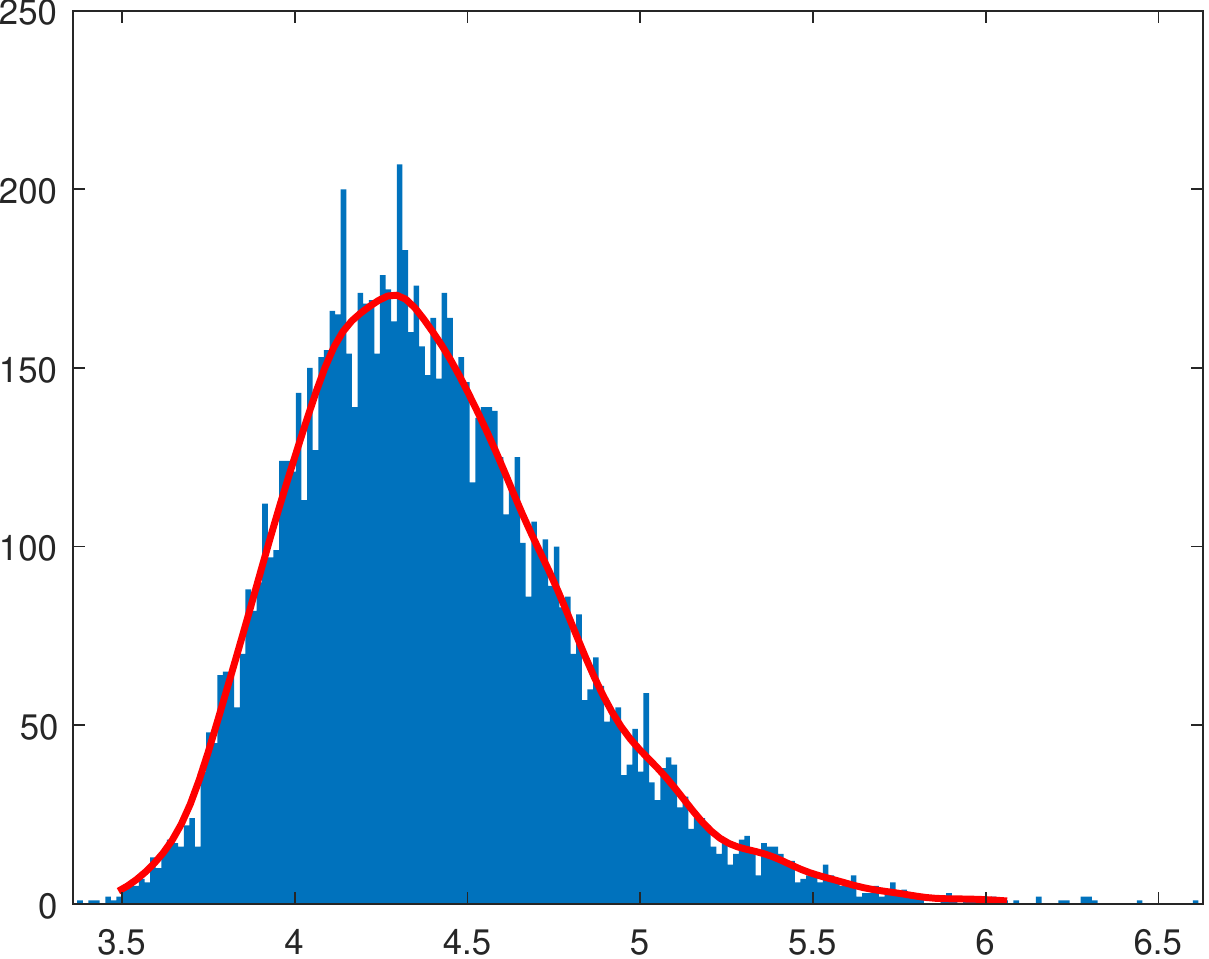}
         \caption{Post-Surgery $\|\tilde{A}^{-1}\|$}
         \label{subfig:Post-Surgery normi}
     \end{subfigure}
     \hspace{0.2cm}
     \begin{subfigure}[b]{0.24\textwidth}
         \centering
         \includegraphics[width=\textwidth]{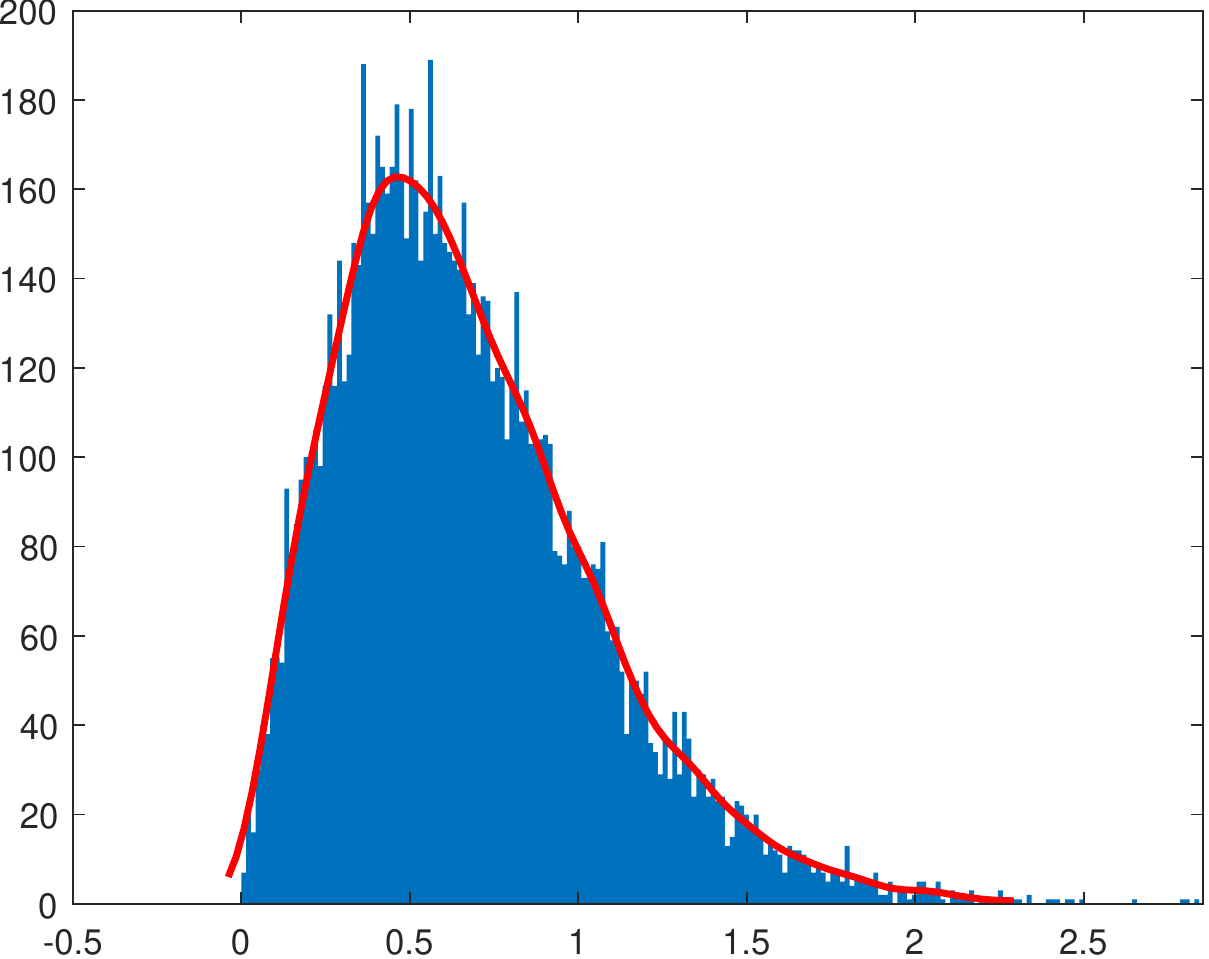}
         \caption{Post-Surgery $\kappa (\tilde{A})$}
         \label{subfig:Post-Surgery CondN}
     \end{subfigure}
        \caption{Distribution of a set of $3 \times 3$ matrices pre and post-surgery: (a) and (d) norm of original matrices, (b) and (e) norm of their inverse, (c) and (f) condition number of matrices}
        \label{fig:Histogram}
\end{figure}

\begin{figure}[h]
     \centering
     \begin{subfigure}[b]{0.33\textwidth}
         \centering
         \includegraphics[width=\textwidth]{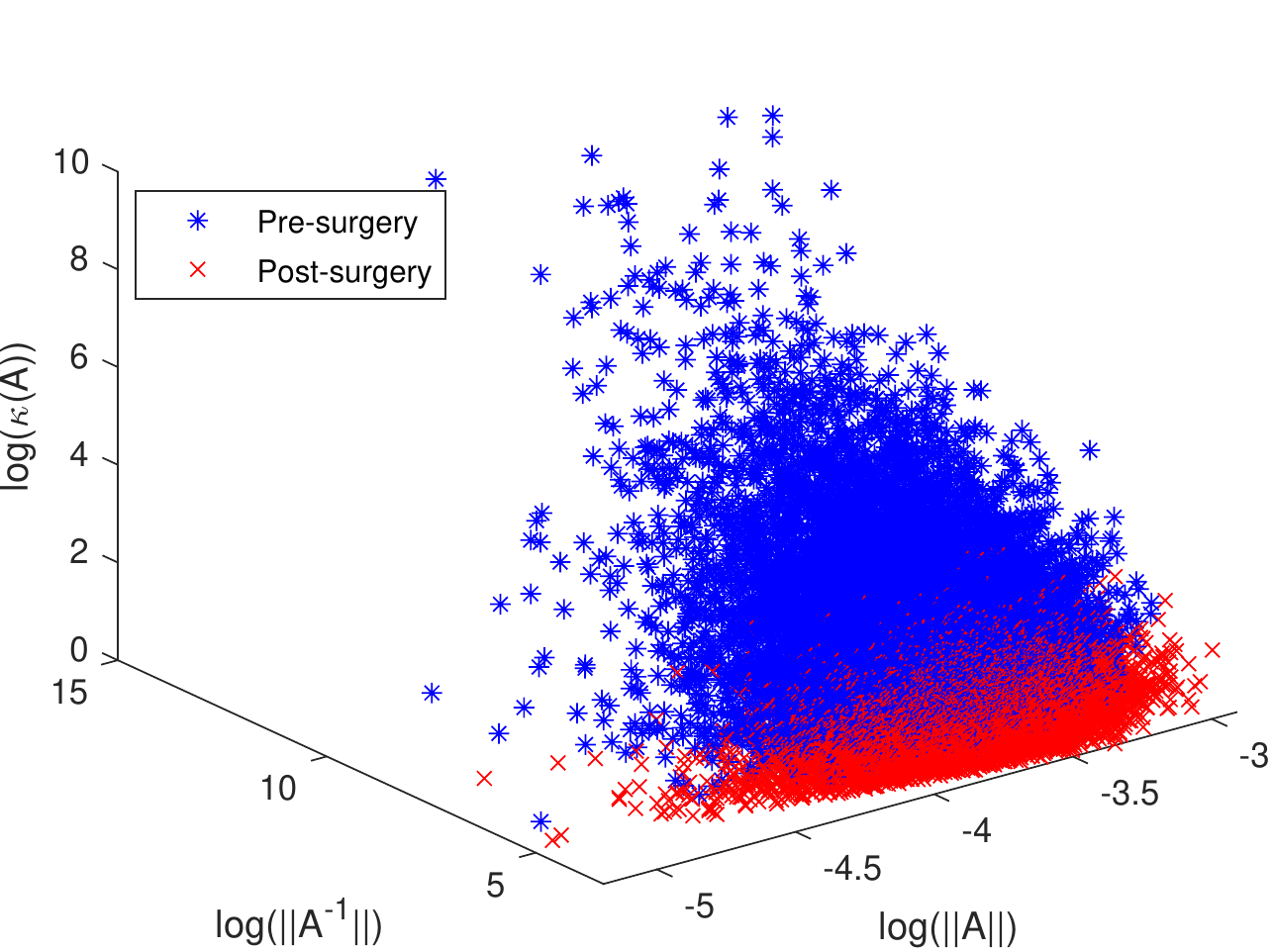}
         \caption{$\tilde{\sigma}_3 $ = $ \sigma_2$}
     \end{subfigure}
     \hfill
     \begin{subfigure}[b]{0.33\textwidth}
         \centering
         \includegraphics[width=\textwidth]{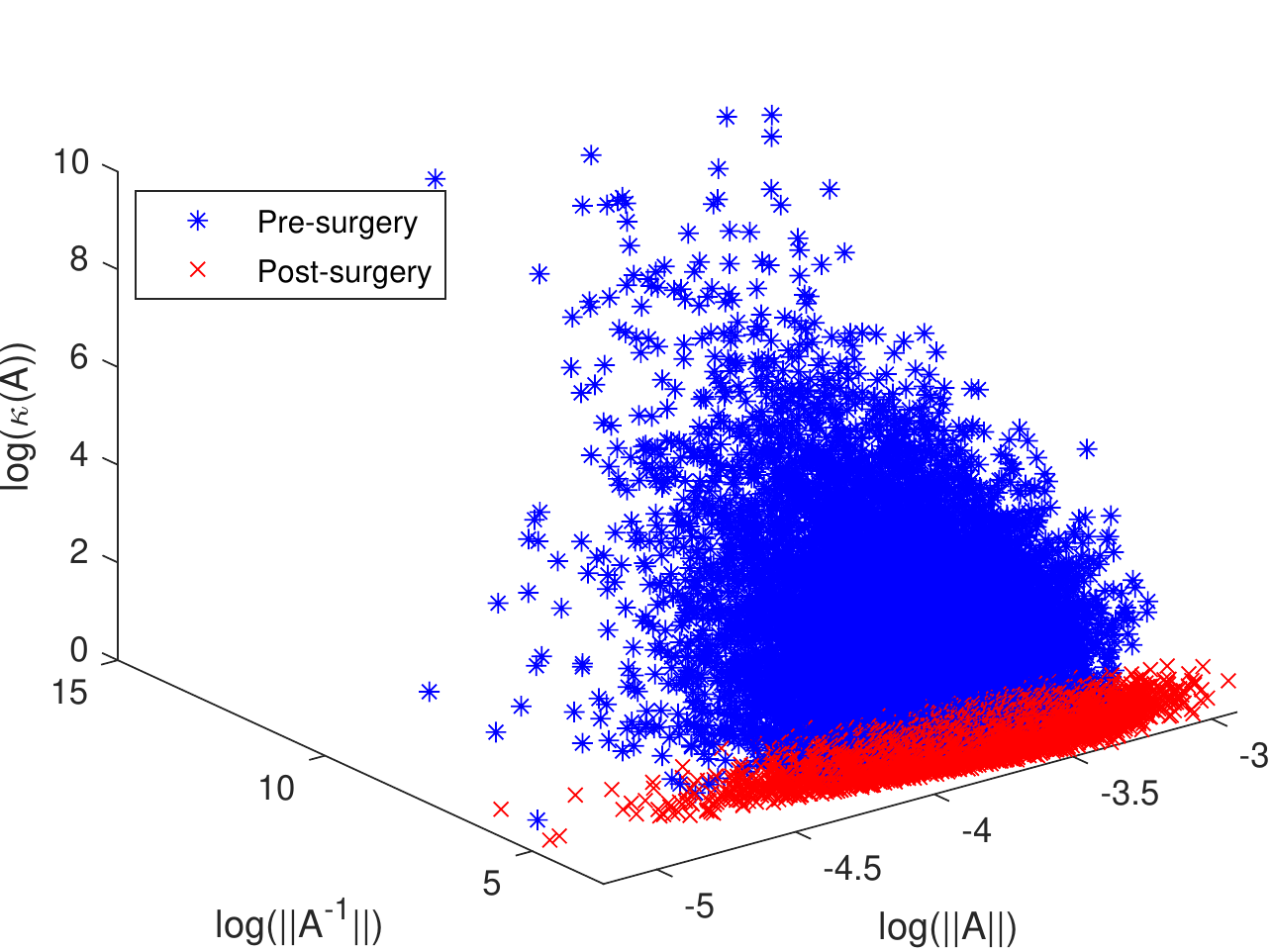}
         \caption{$\tilde{\sigma}_2$ = $\tilde{\sigma}_3$ = $ \sigma_1 /3+ 2\sigma_2 /3$}
     \end{subfigure}
     \hfill
     \begin{subfigure}[b]{0.33\textwidth}
         \centering
         \includegraphics[width=\textwidth]{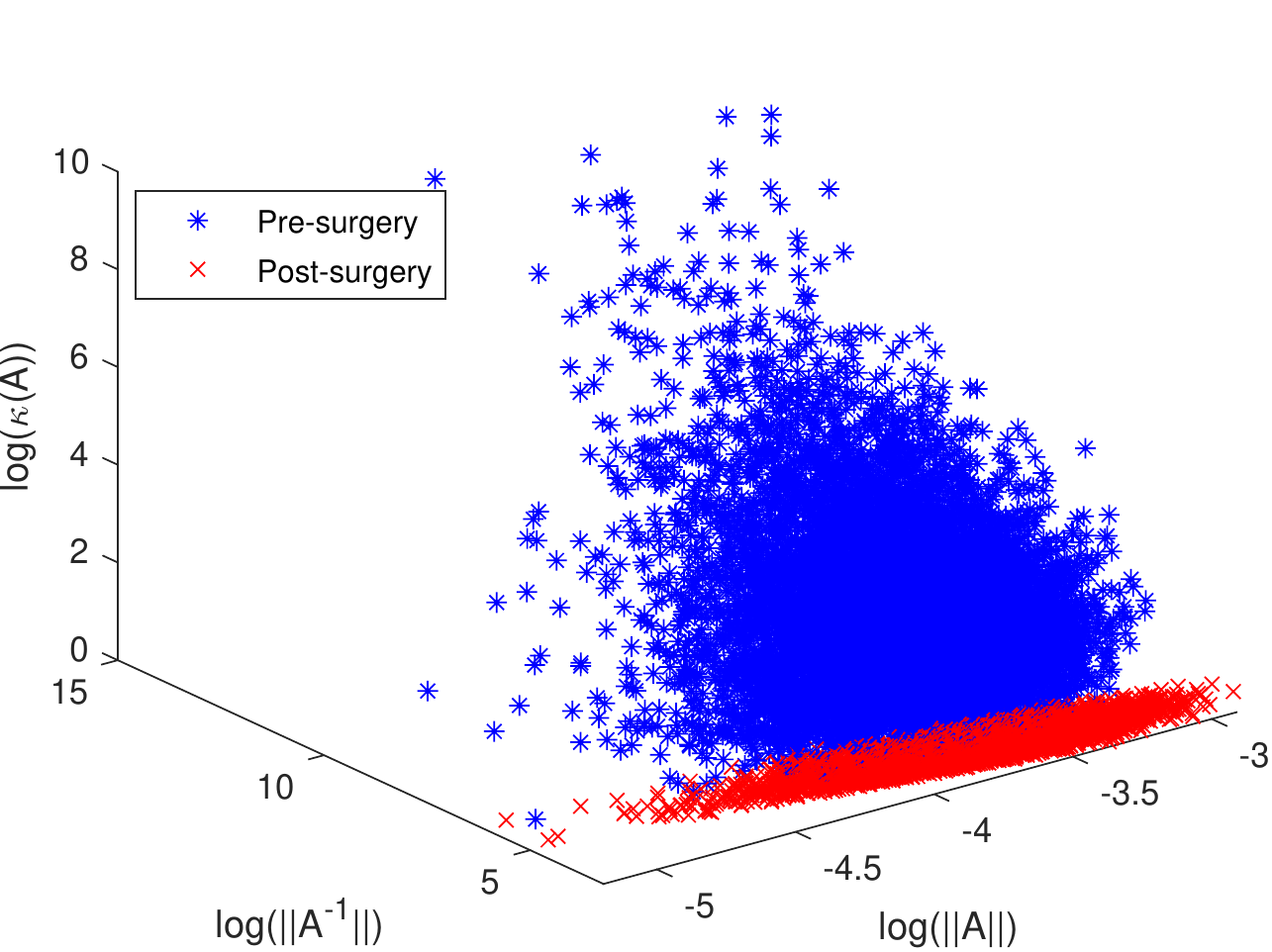}
         \caption{$\tilde{\sigma}_2$ = $\tilde{\sigma}_3 $ = $(\sigma_1 + \sigma_2) /2$}
     \end{subfigure}
     \caption{An illustration of a set of $3 \times 3$ random Gaussian matrices pre- and post-matrix surgery. Shows norm, norm of inverse and condition number (log), where in (a) the smallest singular value $\sigma_3$ is replaced with $\sigma_2$, (b) and (b) $\sigma_2$ and $\sigma_3$ are replaced with new linear combination of $\sigma_1$ and $\sigma_2$.}
     \label{fig:3D plot pre/post-surgery}
\end{figure}

To illustrate the effect of SVD-Surgery on point clouds of convolutions, we generate $10^4$ of $3 \times 3$ matrices drawn from the normal distribution $\mathcal{N}(0,0.01)$. We use the norm of the original matrix, the norm of the inverse, and the condition number to evaluate the change and observe the distribution of these parameters per set. Figure \ref{fig:Histogram} shows a clear difference between the condition number of original versus modified matrices. The reduction in condition number is a result of reducing norms of the inverse of matrices (see figure \ref{fig:3D plot pre/post-surgery}). The minimum and maximum condition numbers for the original set are approximately $1.2$ and $10256$, respectively. After only replacing the smallest singular value $\sigma_3$ with $\sigma_2$ then after reconstructing new minimum and maximum values are 1.006 and 17.14. 

Figure \ref{fig:Histogram} shows a significant change in the distribution of the norm of the inverse of $3 \times 3$ matrices post-surgery and consequently in the distribution of condition numbers. The linear combination is allowing to keep the range of condition number below a certain threshold depending on the distribution of singular values. For instance, 3D illustrations in figure \ref{fig:3D plot pre/post-surgery} show a significant reduction in condition number by keeping the ranges below $3$ in (b) and $2$ in (c), where $\sigma_2$ and $\sigma_3$ are replaced with  $ \sigma_1 /3+ 2\sigma_2 /3$ and $(\sigma_1 + \sigma_2)/2$, respectively. The new minimum and maximum condition number values for both sets after matrix surgery are $[1.004,2.687]$ and $[1.003,1.88]$.

\subsection{Effects of SVD-Surgery on PD’s of convolutional filters point clouds}
One way to control the condition number of CNN filters during training is by controlling their singular values. The implementation of SVD-Surgery can be integrated into customised CNN models for the analysis of natural as well as US image datasets. as a filter regulariser. It can be applied at the filter initialisation when training from scratch, on pretrained filters when training in the transfer learning, as well as on filters modified during training by back-propagation post every batch/epoch. 

In this section, we investigate the topological behaviour of a set of matrices as a point cloud by using the persistent homology tools, as discussed in section \ref{sec:TDA}. For any size $n \times n$ filters, we generate first a set of random Gaussian $n \times n$ matrices and by normalising their entries flattening we get a point cloud in $\mathbb{S}^{n \times n -1}$ residing on its $(n \times n - 1)-sphere$. We construct a second point cloud in $\mathbb{S}^{n \times n - 1}$ by computing the inverse matrices, normalising their entries and flattening. Here, we only illustrate this process for a specific point cloud of $3 \times 3$ matrices for two different linear combinations of the $2$ lower singular values. The general case of larger size filters will be discussed in the first author’s PhD thesis under preparation.

\begin{figure}[h]
     \centering
     \begin{subfigure}[b]{0.3\textwidth}
         \centering
         \includegraphics[width=\textwidth]{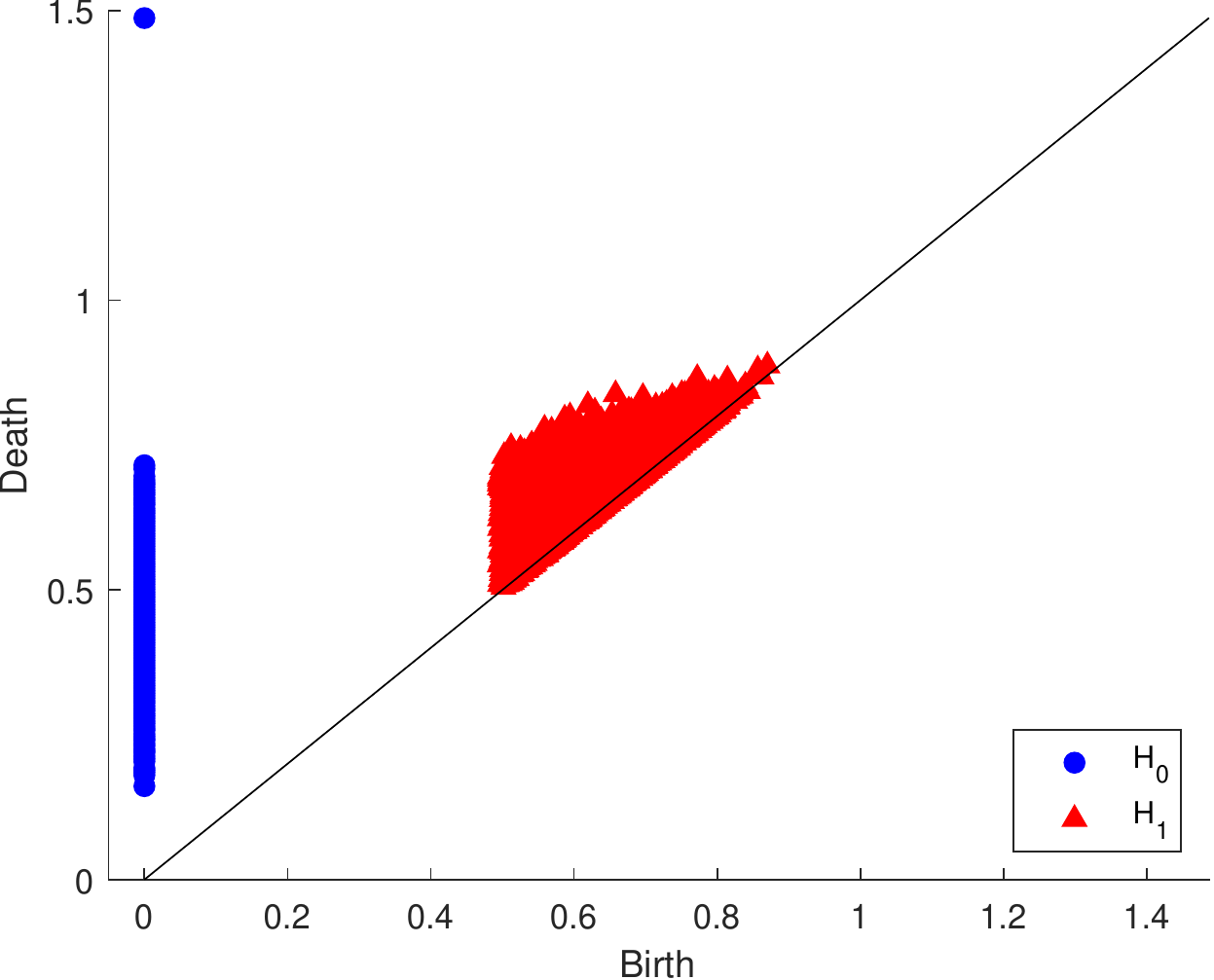}
         \caption{Pre-surgery $\mathcal{A}$}
     \end{subfigure}
      \hfill
     \begin{subfigure}[b]{0.3\textwidth}
         \centering
         \includegraphics[width=\textwidth]{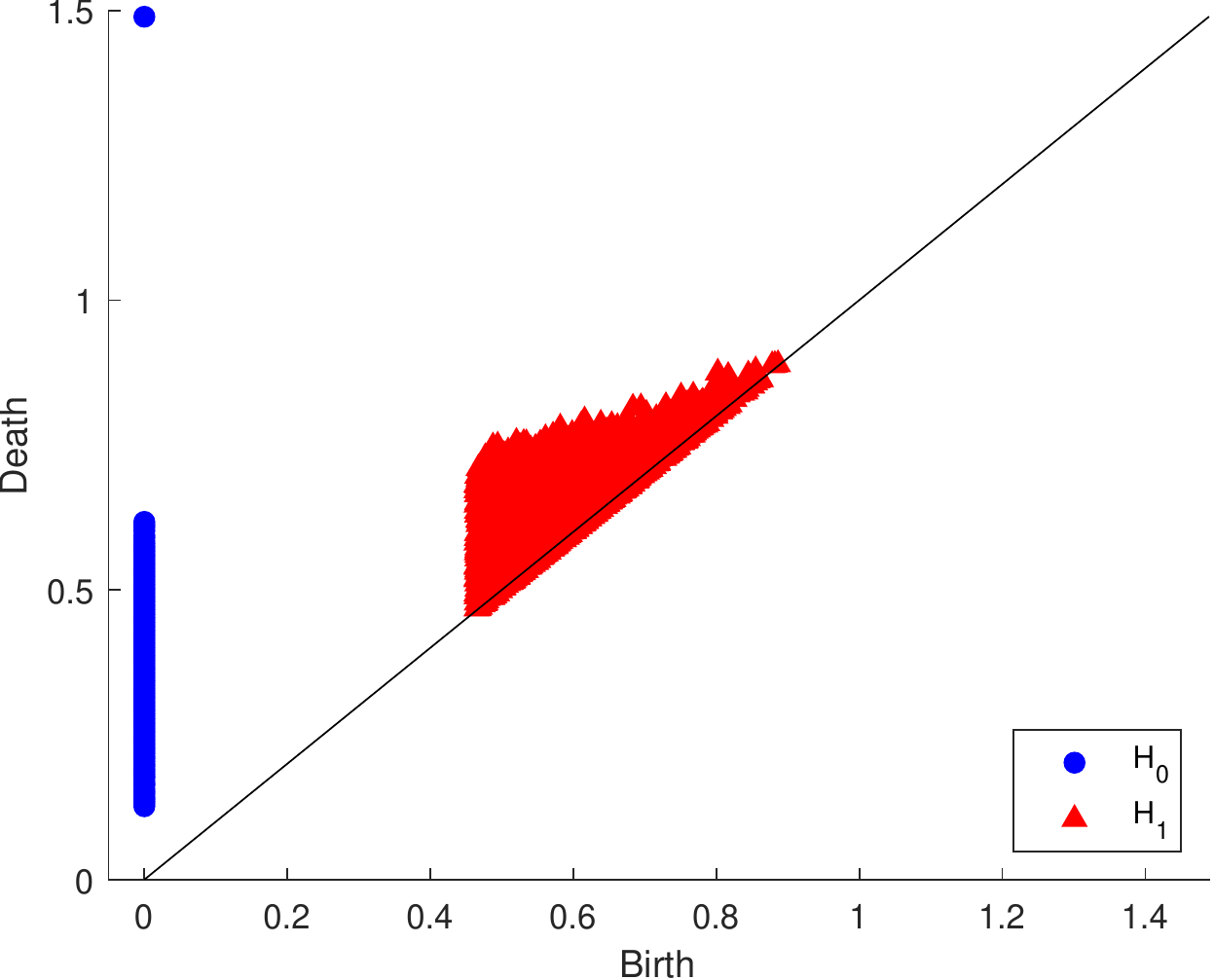}
         \caption{Post-surgery $\tilde{\mathcal{A}}_1$}
     \end{subfigure}
     \hfill
     \begin{subfigure}[b]{0.3\textwidth}
         \centering
         \includegraphics[width=\textwidth]{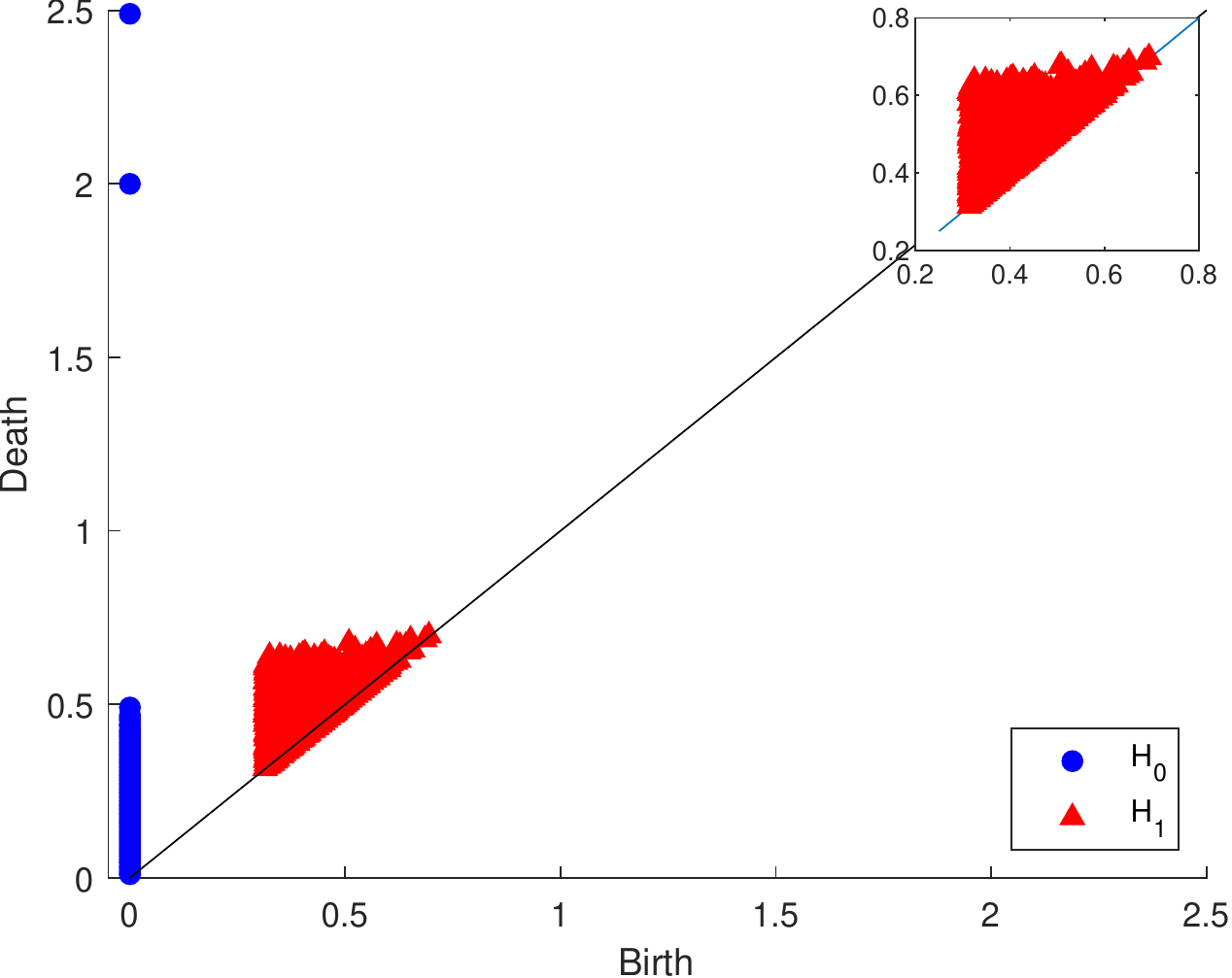}
         \caption{Post-surgery $\tilde{\mathcal{A}}_2$}
     \end{subfigure}
     \hfill
     \begin{subfigure}[b]{0.3\textwidth}
         \centering
         \includegraphics[width=\textwidth]{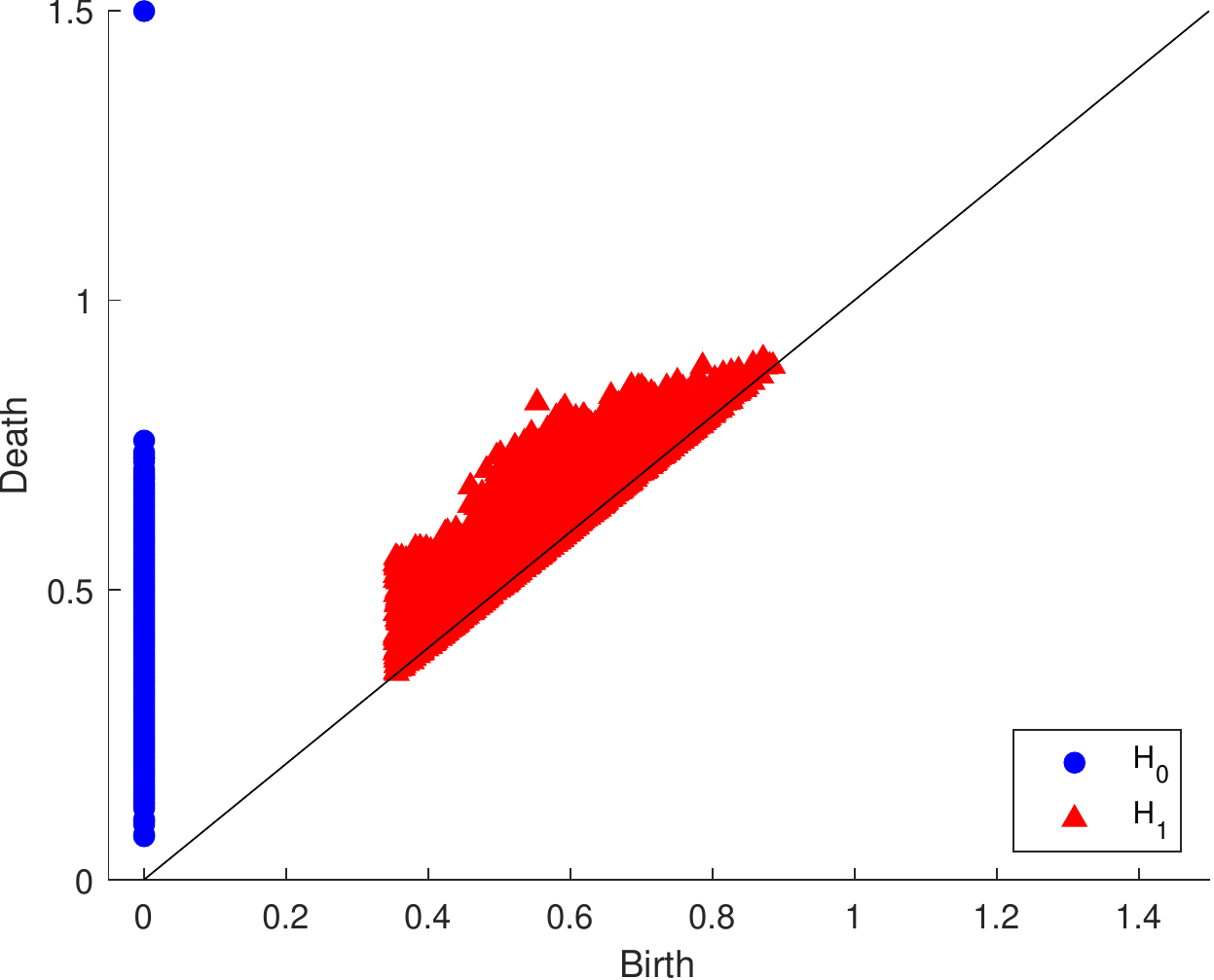}
         \caption{Pre-surgery $\mathcal{A}^{-1}$}
     \end{subfigure}
     \hfill
     \begin{subfigure}[b]{0.3\textwidth}
         \centering
         \includegraphics[width=\textwidth]{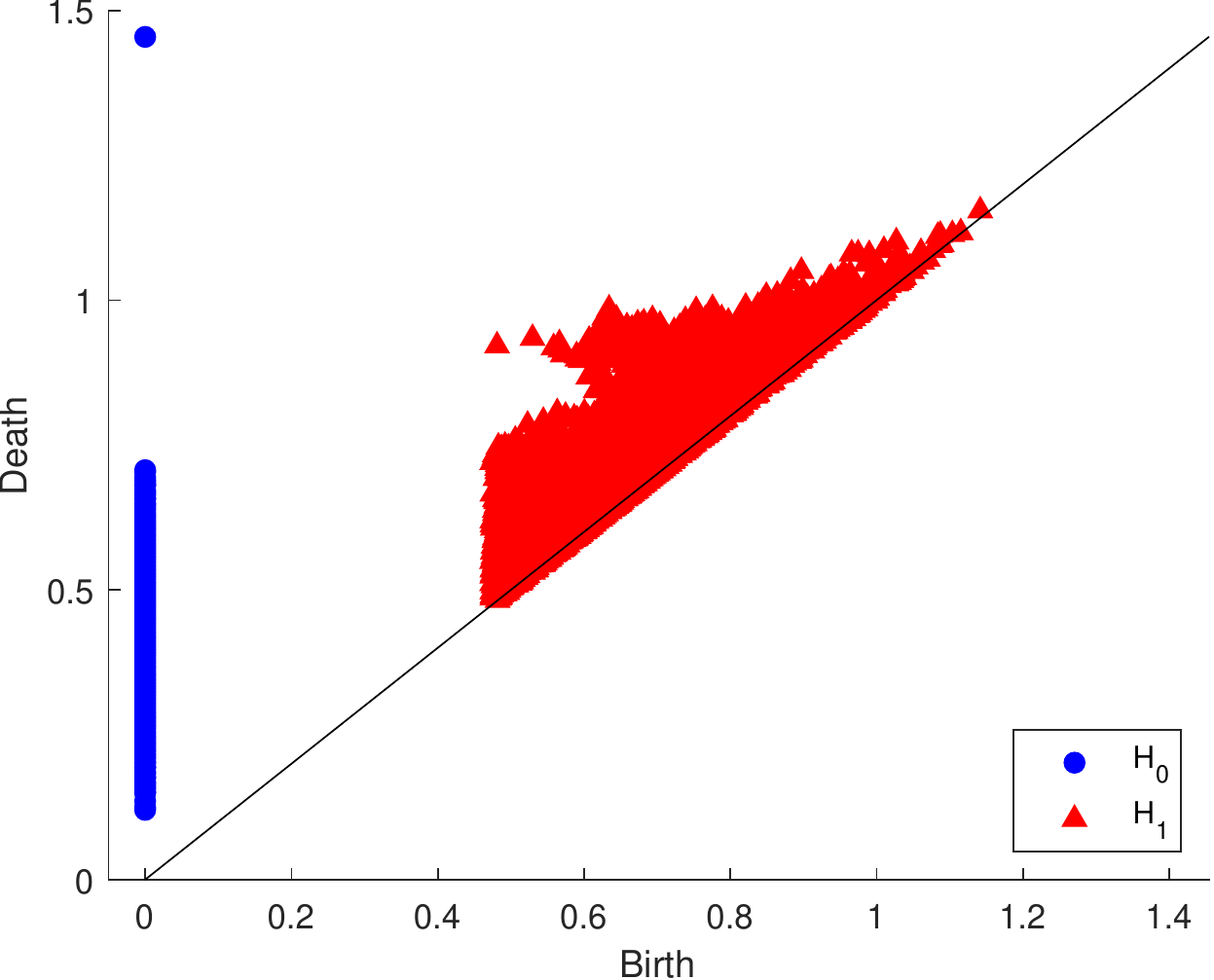}
         \caption{Post-surgery $\tilde{\mathcal{A}}^{-1}_1$}
     \end{subfigure}
     \hfill
     \begin{subfigure}[b]{0.3\textwidth}
         \centering
         \includegraphics[width=\textwidth]{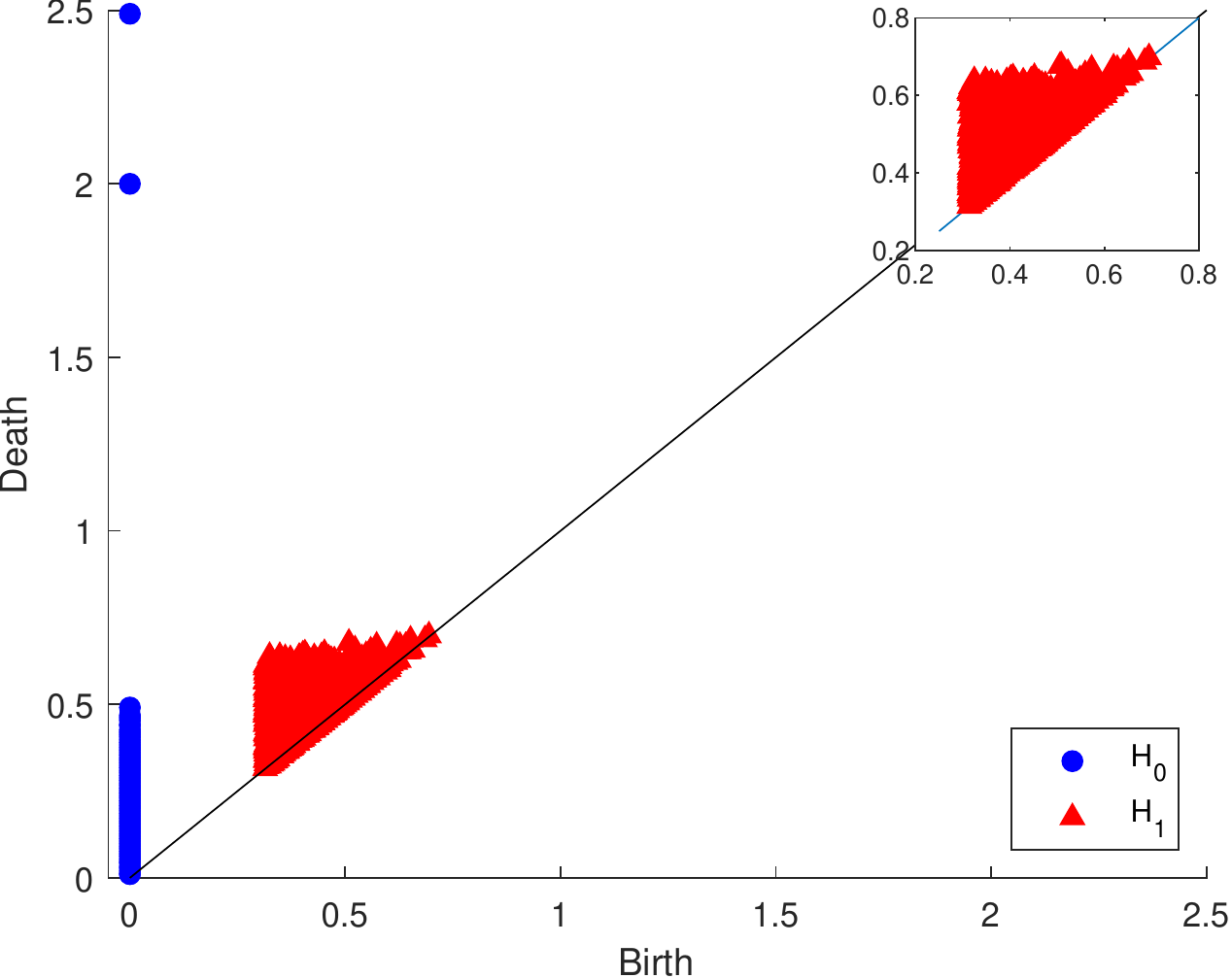}
         \caption{Post-surgery $\tilde{\mathcal{A}}^{-1}_2$}
     \end{subfigure}
        \caption{Persistence diagram of point clouds $\mathcal{A}$ and $\mathcal{A}^{-1}$ before and after SVD based surgery.}
        \label{fig:PD pre and post-surgery}
\end{figure}

Figure \ref{fig:PD pre and post-surgery}, below, shows the $H_0$ and $H_1$, persistence diagrams for point clouds (originals and inverses) plus those for post-matrix surgery with respect to the linear combinations: (1) replacing both $\sigma_2$ and $\sigma_3$ with $\sigma_1$  (i.e $\kappa (A)=1$), and (2) replacing $\sigma_3$ with $\sigma_1$. The first row, corresponds to the effect of SVD on the PD of the original point cloud, while the second row corresponds to the inverse point cloud.

It is clear that PDs of $\tilde{\mathcal{A}}_2$ and $\tilde{\mathcal{A}}_2^{-1}$ are equivalent as a reflection of the fact that this surgery produces optimally well-conditioning matrices that are orthogonal. Note that in this case, the inverse matrices are simply the transpose of the original ones. Point clouds $\tilde{\mathcal{A}}_1$ and $\tilde{\mathcal{A}}_1^{-1}$ are consist of a set of matrices $\tilde{A}_1$ and $\tilde{A}_1^{-1}$ where the smallest singular value $\sigma_3$ is replaced with $\sigma_2$ and the range of condition number from $[1.2,10256]$ to $[1.006, 17.14]$. It is customary, to determine differences and similarities of PDs using distance measures such as bottleneck distance.

Our ongoing investigations reveal noticeable improvements in CNN model performance regarding robustness to tolerable data perturbations, and generalisation to unseen data when using SVD-Surgery for convolutional layer filters at (1) initialisation from scratch, (2) pretrained filters, (3) during training batches and/or epochs, and particularly (4) when all combined. Furthermore, controlling condition numbers during training stabilises topological behaviour of filters per convolutional layer. Future publications, will be dealing with these applications of matrix surgery and topological data analysis in CNNs.

\section{Conclusion}
We introduced a simple strategy for matrix surgery to reduce and control the condition number of $n \times n$ matrices, by replacing all the singular values of the input matrix starting with a selected diagonal entry downward with a convex linear combination of the list, and reconstructing the matrix with a much lower condition number. We demonstrated that this strategy applied on several point clouds of sufficiently large convolution filters preserve filters’ norm and reduces the norm of its inverse depending on the chosen linear combination parameters. Our approach showed significant improvements in condition number and persistent homology. In terms of our motivating challenge, our SVD-Surgery approach is ideally efficient in controlling condition numbers of convolutional layer filters during training (from scratch or in transfer learning mode) without increasing model complexity or learning extra hyper-parameters. It preserves monotonicity of the singular values and the linear combination parameters can be made layer dependent determined to avoid making significant rescaling of training dataset features along the singular vectors.   

\section*{Acknowledgement}
This research is sponsored by Ten-D AI Medical Technologies Ltd.

\bibliographystyle{ieeetr}  
\bibliography{references}  
\end{document}